\documentclass[10pt, final, oneside, twocolumn, letter]{IEEEtran}

\usepackage{graphicx}           
\usepackage{mathrsfs}
\usepackage{amssymb}
\usepackage{amsmath,mathtools}
\usepackage{amsmath}
\usepackage{amsthm}
\usepackage{tikz}
\usepackage{caption}
\usepackage{subfig}
\usepackage{enumerate}
\usepackage{epstopdf}
\usepackage{array}
\newcolumntype{P}[1]{>{\raggedleft\arraybackslash}p{#1}}
\everymath{\displaystyle}

\newcommand{\lls}{\langle\langle}
\newcommand{\ggs}{\rangle\rangle}

\newtheorem{theorem}{\textbf{Theorem}}
\newtheorem*{theorem*}{\textbf{Theorem}}
\newtheorem{definition}{\textbf{Definition}}
\newtheorem{prop}{\textbf{Proposition}}
\newtheorem{assum}{\textbf{Assumption}}

\newtheorem{lemma}{Lemma}

\title{Construction of Synergistic Potential Functions on SO(3) with Application to Velocity-Free Hybrid Attitude Stabilization}
\author{Soulaimane Berkane and Abdelhamid Tayebi
\thanks{This work was supported by the National Sciences and Engineering Research Council of Canada (NSERC).}
\thanks{The authors are with the Department of Electrical and Computer Engineering, University of Western Ontario, London, Ontario, Canada. The second author is also with the Department of Electrical Engineering, Lakehead University, Thunder Bay, Ontario, Canada.
     {\tt\small sberkane@uwo.ca, tayebi@ieee.org} }%
}
\begin{document}
\maketitle
\begin{abstract}
We propose a systematic and comprehensive procedure for the construction of synergistic potential functions, which are instrumental in hybrid control design on $SO(3)$. A new map via angular warping on $SO(3)$ is introduced for the construction of such a family of potential functions allowing an explicit determination of the critical points and the synergistic gap. Some optimization results on the synergistic gap are also provided. The proposed synergistic potential functions are used for the design of a global velocity-free hybrid attitude stabilization scheme relying solely on inertial vector measurements. Comparative simulation results between the proposed global hybrid control scheme and the almost global smooth control scheme have been carried out.

\end{abstract}
\section{Introduction}
The attitude control problem of rigid body systems has been widely treated in the literature with many applications related to aerospace and marine engineering  (see, for instance, \cite{Meyer1971}, \cite{Tayebi2006}, \cite{Tayebi2008a} and \cite{Mayhew2011b}). The major challenge of this control problem is related to the motion space $SO(3)$ where the angular velocity is not a straightforward derivative of the angular position. In fact, it was shown in \cite{Bhat2000} that, due to the inherent topology of the compact manifold $SO(3)$, there is no continuous time -invariant feedback that globally stabilizes the rigid body attitude to a desired reference. The best results that one can achieve on $SO(3)$, or any Lie group diffeomorphic to $SO(3)\times\mathbb{R}^3$, with time-invariant smooth feedback laws are \textit{almost} global, where the attitude is stabilized from any initial attitude except from a set of Lebesgue measure zero (see, for instance, \cite{Koditschek}, \cite{Sanyal2009}, \cite{Chaturvedi2011}). Most attitude control systems are obtained using the popular modified-trace function
$
V_A(R)=\mathrm{tr}(A(I-R)),
$
where $\mathrm{tr}(A)I-A$ is symmetric positive definite, in the Lyapunov design and analysis. For instance, consider the attitude kinematics
$
\dot R=R\Omega,
$
then the angular velocity tensor $\Omega\in\mathfrak{so}(3)$ can not be designed to globally stabilize the attitude $R\in SO(3)$, at any given reference, if it is restricted to be smooth.
Nevertheless, under the continuous feedback
\begin{equation*}
\Omega=-R^T\nabla V_A(R)=-(AR-R^\top A),
\end{equation*}
one can show that the closed loop system has multiple equilibria corresponding to the critical points of $V_A$ where the gradient $\nabla V_A(R)$ is zero. The \textit{undesired} critical points are identified by $\mathcal{R}_a(\pi,\mathcal{E}(A))$ which represents the set of all rotations of angle $\pi$ and axis $v\in\mathcal{E}(A)$, with $\mathcal{E}(A)$ being the set of unit eigenvectors of $A$. It is also shown that the manifold $\mathcal{R}_a(\pi,\mathbb{S}^2)$ of all rotations of angle $\pi$ is invariant under this feedback \cite{mayhew2011synergistic}. The appearance of undesired critical points when considering smooth controls on $SO(3)$ is non-avoidable. This is mainly due to the fact that, according to Morse theory \cite{morse1934calculus}, any smooth potential function on $SO(3)$ is guaranteed to have at least four critical points where its gradient vanishes. 

On the other hand, there have been some attempts to design attitude control systems with global stability results by introducing discontinuities. For instance, using a discontinuous quaternion-based control, as done in \cite{Thienel2003}, one can achieve global stability results. However, these discontinuous attitude control systems, in addition to the quaternion ambiguity, suffer from non-robustness to arbitrary small measurement disturbances as explained in \cite{Mayhew2011a}.

The recent work in \cite{Mayhew2011} focuses on the design of hybrid feedback systems that are able to overcome the topological obstruction to global asymptotic stability on $SO(3)$ while, in the same time, ensuring some robustness to measurement noise. The hybrid algorithm is based on a family of potential functions and a hysteresis-based switching mechanism to avoid the undesired critical points. After each switching, the control law derived from the minimal potential function is selected. A sufficient condition for global asymptotic stability of the resulting hybrid controller is the ``synergism" property of this family of potential functions. A family of potential functions on $SO(3)$ is synergistic if at each critical point (other than the desired one) of a potential function in the family, there exists another potential function in the family that has a lower value. Moreover, if all the potential functions in the family
share the identity element $I_{3\times 3}$ as a critical point then it is called a \textit{centrally} synergistic family. Thanks to the hysteresis gap, this controller guarantees robustness to small measurements noise.

As a consequence of this approach, the design of hybrid controllers on $SO(3)$, leading to robust and global asymptotic stability results, boils down to the search for a \textit{suitable} synergistic family of potential functions. The work in \cite{mayhew2011synergistic} suggests the technique of ``angular warping" to construct synergistic potential functions on $SO(3)$, although without rigorous proof of synergism. In \cite{casau2014globally}, the authors proved that this technique could generate a synergistic family, under some conditions, when applied to the modified trace function. The major drawback of these approaches is related to the difficulty of determining the synergistic gap which is required for the implementation of the hybrid controller. In an attempt to solve this problem, the authors in \cite{Mayhew2013} tried to relax the centrality assumption by considering scaled, biased and translated modified trace functions. However, the sufficient synergism conditions provided therein were very conservative, difficult to satisfy and only hand tuning of the parameters was proposed. Another form of non-central synergistic potential functions appeared in \cite{lee2015tracking} by comparing the actual and desired directions, leading to a simple expression of the synergistic gap.

In this paper, we consider a central family of potential functions obtained via angular warping on $SO(3)$ and derive necessary and sufficient conditions guaranteeing that the family under consideration is synergistic. We propose a new warping angle function that allows explicit calculation of the critical points as well as the synergistic gap. We also provide sufficient conditions on the angular warping direction to maximize the synergistic gap. The fact that our approach generates a \textit{central} synergistic family is advantageous. In fact, each control law, derived from each smooth potential function in the central
family, guarantees (independently) almost global asymptotic stabilization of the attitude. This is desirable in practice since the control objective remains achievable even when the hybrid switching mechanism runs into error. Moreover, the control law derived from the presented central synergistic family is directly expressed in terms of an \textit{arbitrary} number of vector measurements and \textit{only} two potential functions are sufficient to guarantee synergism. It is also worth pointing out that our approach to the construction of synergistic families on $SO(3)$ can be adapted to other compact manifolds such as $\mathbb{S}^1$ and $\mathbb{S}^2$ where angular warping technique has been used \cite{mayhew2010planar,mayhew2010spherical, mayhew2010pendulum}. 
\\
The second contribution consists in the design of a hybrid velocity-free attitude stabilization scheme relying solely on inertial vector measurements. The proposed synergistic potential functions were instrumental in the design of such a hybrid controller leading to global asymptotic stability results. The proposed control scheme is inspired from our earlier work \cite{TAC2013} where \textit{almost} global asymptotic stability results have been obtained. Note that in \cite{Mayhew2011b} a hybrid velocity-free attitude controller, inspired by \cite{Tayebi2008a}, has been proposed assuming that the attitude is available for feedback. Since there is no sensor that provides directly a measurement of the attitude, an attitude estimation algorithm, which usually relies on angular velocity measurements, is required.
\section{Attitude reprentation and Mathematical Preliminaries}\label{Att_rep}
\subsection{Notations and mathematical preliminaries}\label{notation}
The sets of real, non-negative real and natural numbers are denoted as $\mathbb{R}$, $\mathbb{R}^+$ and $\mathbb{N}$, respectively. $\mathbb{R}^n$ denotes the $n$-dimensional Euclidean space and $\mathbb{S}^n$ denotes the unit $n$-sphere embedded in $\mathbb{R}^{n+1}$. Given two matrices $A,B\in\mathbb{R}^{m\times n}$, their Euclidean inner product is defined as $\lls A,B\ggs=\textrm{tr}(A^{\top}B)$ where $(\cdot)^{\top}$ denotes the transpose of $(\cdot)$. The $2$-norm of a vector $x\in\mathbb{R}^n$ is $\|x\|=\sqrt{x^\top x}$ and the Frobenius norm of a matrix $A\in\mathbb{R}^{n\times m}$ is $\|A\|_F=\sqrt{\lls A,A\ggs}$. For a given square matrix $A\in\mathbb{R}^{n\times n}$, we denote by $\lambda_i^A, \lambda_{\mathrm{min}}^A\;$ and $\lambda_{\mathrm{max}}^A$ the $i$-th, minimum and maximum eigenvalue of $A$, respectively.
\\
Given a manifold $M$, a tangent vector at $x\in M$ is $\gamma'(0):=d\gamma(\tau)/d\tau|_{\tau=0}$ for some smooth path $\gamma: \mathbb{R}\to M$ such that $\gamma(0)=x$. The \textit{tangent space} to $M$ at $x$ is the set of all tangent vectors at $x$, denoted $T_xM$. The disjoint union of all tangent spaces forms the \textit{tangent bundle} $TM$.
Let $M$ and $N$ be two smooth manifolds and let $f: M\to N $ be a differentiable map. The \textit{tangent map (differential)} of $f$ at a point $x\in M$ is the map \cite{Oxford2009}
$$\begin{array}{rl}
df(x): T_xM&\to T_{f(x)}N\\
				\xi&\mapsto df(x)\cdot\xi:=\left(f\circ\gamma\right)^\prime(0),
\end{array}
$$
where $\gamma(\tau)$ is a path in $M$ such that $\gamma(0)=x$ and $\gamma^\prime(0)=\xi$. The inverse image of a subset $\mathcal{S}_N\subseteq N$ under the map $f$ is the subset of $M$ defined by
$
f^{-1}(\mathcal{S}_N)=\{x\in M\mid f(x)\in\mathcal{S}_N\}.
$
Let $f: M\to\mathbb{R}$ be a differentiable real-valued function. A point $x\in M$ is called a \textit{critical point}\footnote{For a reference, see Morse Theory 279 (VII.16), page 1049 of \cite{ito1993encyclopedic}.} of $f$ if the differential map $df(x)\cdot\xi$ is zero at $x$ for all $\xi\in T_xM$. We denote by $\Psi_f\subseteq M$ the set of all critical points of $f$ on $M$. Let $\langle\;,\;\rangle_x: T_xM\times T_xM\to \mathbb{R}$ be a \textit{Riemannian metric} on $M$. The \textit{gradient} of $f$, denoted $\nabla f(x)\in T_xM$, relative to the Riemannian metric $\langle\;,\;\rangle_x$ is uniquely defined by
\begin{equation}\label{grad_def}
df(x)\cdot\xi=\langle \nabla f(x),\xi\rangle_x\;\;\;\mathrm{for\;all\;}\xi\in T_xM.
\end{equation}
\subsection{Attitude representation and kinematics}\label{Sec_attitude}
Consider the \textit{general linear group} $GL(3)$ which represents the set of $3\times 3$ invertible matrices, together with the ordinary matrix multiplication. A square matrix $R\in GL(3)$ is called a \textit{rotation matrix} if $R$ belongs to the \textit{special orthogonal group} $SO(3)\subset GL(3)$ where
$
SO(3) := \{ R \in \mathbb{R}^{3\times 3}|\mathrm{det}(R)=1,RR^{\top}=I\},
$
and $I $ is the three-dimensional identity matrix. The \textit{Lie algebra} of $SO(3)$, denoted by $\mathfrak{so}(3)$,  is the vector space of 3-by-3 skew-symmetric matrices
$
\mathfrak{so}(3)=\left\{\Omega\in\mathbb{R}^{3\times 3}\mid\;\Omega^{\top}=-\Omega\right\}.
$
The group $SO(3)$ has a compact manifold structure where its \textit{tangent spaces} are identified by
$
T_RSO(3):=\left\{R\Omega\mid\Omega\in\mathfrak{so}(3)\right\}.
$
The Euclidean inner product on $\mathbb{R}^{3\times 3}$, when restricted to the Lie-algebra of skew symmetric matrices, defines the following \textit{left-invariant} Riemannian metric on $SO(3)$
\begin{equation}\label{metric}
\langle R\Omega_1,R\Omega_2\rangle_R:=\lls\Omega_1,\Omega_2\ggs, \;\;\;
\end{equation}
for all $R\in SO(3)$ and $\Omega_1,\Omega_2\in\mathfrak{so}(3).$
Let $\times$ denotes the vector cross-product on $\mathbb{R}^3$ and define the map $[~.~]_\times: \mathbb{R}^3\to\mathfrak{so}(3)$; $\omega \mapsto [\omega]_\times$ such that
$
[\omega]_\times u=\omega\times u,$ for all $\omega, u\in\mathbb{R}^3,
$
where for any vector $\omega\in\mathbb{R}^3$, we have
\begin{equation*}
[\omega]_\times= \left[\begin{array}{ccc}
0 & -\omega_3 & \omega_2 \\
\omega_3 & 0 &-\omega_1 \\
-\omega_2 &\omega_1 &0
\end{array}\right].
\label{skew}
\end{equation*}
If $a=[\alpha_1, \alpha_2, \alpha_3]^{\top}$ and $b=[\beta_1, \beta_2, \beta_3]^{\top}$ are vectors in $\mathbb{R}^3$, represented in some orthonormal basis $\mathcal{B}=\left\{v_1, v_2, v_3\right\}$, then their cross product can be written as \cite{arfken1999mathematical}
\begin{equation}\label{a_cross_b}
a\times b=\sum_{m,n,l}\varepsilon_{mnl}\alpha_m\beta_nv_l,
\end{equation}
where $\varepsilon_{mnl}$ is the Levi-Cevita symbol defined by
\begin{equation*}
\varepsilon_{mnl}=\left\{\begin{array}{l}
0\;\;\;\mathrm{for\;}m=n, m=l \mathrm{\;or\;}n=l\\
+1\;\;\;\mathrm{for\;}(m,n,l)\in\left\{(1,2,3), (2,3,1), (3,1,2)\right\}\\
-1\;\;\;\mathrm{for\;}(m,n,l)\in\left\{(1,3,2), (3,2,1), (2,1,3)\right\}.
\end{array}\right.
\end{equation*}
Let $\mathrm{vex}:\mathfrak{so}(3)\to\mathbb{R}^3$ denotes the inverse isomorphism of the map $[~.~]_\times$, such that
$
\mathrm{vex}([\omega]_\times)=\omega,$ for all $\omega\in\mathbb{R}^3$ and $[\mathrm{vex}(\Omega)]_\times=\Omega,$ for all $\Omega\in\mathfrak{so}(3)$. By defining $\mathbb{P}_a:\mathbb{R}^{3\times 3}\to\mathfrak{so}(3)$ as the projection map on the Lie algebra $\mathfrak{so}(3)$ such that
$
\mathbb{P}_a(A):=(A-A^{\top})/2,
$
one can extend the definition of $\mathrm{vex}$ to $\mathbb{R}^{3\times 3}$ by taking the composition map $\psi := \mathrm{vex}\circ \mathbb{P}_a$ 
        such that, for a $3$-by-$3$ matrix $A:=[a_{ij}]_{i,j = 1,2,3}$, one has
        \begin{equation}\label{psi}
        \psi(A):=\mathrm{vex}\left(\mathbb{P}_a(A)\right)=\frac{1}{2}\left[\begin{array}{c}
        a_{32}-a_{23}\\a_{13}-a_{31}\\a_{21}-a_{12}
        \end{array}
        \right].
        \end{equation}
The following identity is used throughout this paper:
\begin{equation}
\label{id7}
\lls A,[u]_\times\ggs=2\psi(A)^{\top}u,
\end{equation}
where $A\in\mathbb{R}^{3\times 3}$ and $u\in\mathbb{R}^3$. Alternatively, an element $R\in SO(3)$ can be represented as a rotation of angle $\theta\in\mathbb{R}$ around a unit vector axis $u\in\mathbb{S}^2$. This is commonly known as the angle-axis parametrization of $SO(3)$ and it is given by the map $\mathcal{R}_a:\mathbb{R}\times\mathbb{S}^2\to SO(3)$:
\begin{equation}\label{Rodrigues}
\mathcal{R}_a(\theta,u):=e^{\theta[u]_\times}=I+\sin(\theta)[u]_\times+(1-\cos(\theta))[u]_\times^2,
\end{equation}
where $e^A$ denotes the matrix exponential of $A$. Equation (\ref{Rodrigues}) is known as Rodriguez formula. In addition, the map $\mathcal{R}_a$ satisfies the following property
\begin{equation}\label{Rt1Rt2}
\mathcal{R}_a(\theta_1,u)\mathcal{R}_a(\theta_2,u)=\mathcal{R}_a(\theta_1+\theta_2,u),\;\;\;\forall\theta_1.\theta_2\in\mathbb{R},\;\forall u\in\mathbb{S}^2.
\end{equation}
In this work, we may also use the unit quaternion representation\footnote{For more details on the unit-quaternion (in addition to other forms of attitude representation) the reader is referred to \cite{Shuster1993}, \cite{murray1994mir}, and \cite{hughes1986sad}.} of a rotation matrix $R\in\ SO(3)$. A unit quaternion $Q = (\eta,\epsilon)\in\mathbb{Q}$, consists of a scalar part $\eta$ and three-dimensional vector $\epsilon$, such that
$
\mathbb{Q}:=\{Q = (\eta,\epsilon)\in\mathbb{R}^4\;|\; \eta^2+\epsilon^{\top}\epsilon=1\}.
$
The relation between the quaternion representation and the angle-axis representation is given by $
\eta=\cos\left(\theta/2\right)$ and $\epsilon=\sin\left(\theta/2\right)u$.
Therefore, a unit quaternion represents a rotation matrix through the map $\mathcal{R}_Q:\mathbb{Q}\to SO(3)$ defined as
\begin{equation}\label{rod}
\mathcal{R}_Q(Q)=I+2[\epsilon]^2_\times+2\eta[\epsilon]_\times.
\end{equation}
The set $\mathbb{Q}$ forms a group with the quaternion product, denoted by $\odot$, being the group operation and quaternion inverse defined by $Q^{-1} = \left(\eta, -\epsilon \right)$ as well as the identity-quaternion $Q = \left( 1, 0_{3\times 1} \right)$, where $0_{3\times 1}\in \mathbb{R}^3$ is a column vector of  zeros. Given $Q_1,Q_2\in \mathbb{Q}$ where $Q_1= (\eta_1,\epsilon_1)$ and $Q_2 = (\eta_2,\epsilon_2)$ the quaternion product is defined by
\begin{equation}\label{Qmultiply}
Q_1\odot Q_2 =\left( \eta_1 \eta_2- \epsilon_1^{\top}\epsilon_2 ,\eta_1\epsilon_2 + \eta_2\epsilon_1 +
[\epsilon_1]_\times\epsilon_2 \right),
\end{equation}
and satisfying \begin{equation}\label{R1R2}
\mathcal{R}_Q(Q_1)\mathcal{R}_Q(Q_2)=\mathcal{R}_Q(Q_1\odot Q_2).
\end{equation}
\subsection{Hybrid Systems Framework}
 In this paper, we make use of the recent framework for dynamical hybrid systems found in \cite{Goebel2006, Goebel2009}. A subset $E\subset\mathbb{R}_{\geq 0}\times\mathbb{N}$ is a  \textit{hybrid time domain}, 
        if it is a union of finitely or infinitely many intervals of the form $[t_j ,t_{j+1}]\times\{j\}$ where $0=t_0\leq t_1\leq t_2\leq...$,  with the last interval being possibly of the form $[t_j ,t_{j+1}]\times\{j\}$ or $[t_j ,\infty)\times\{j\}$. Let $\rightrightarrows$ denote a set-valued mapping. A general model of a hybrid system $\mathcal{H}$ takes the form:
            \begin{equation}\label{Hybrid:general}
            \mathcal{H}\left\{\begin{array}{l}
            \hspace{.25cm}\dot x\in F(x),\hspace{.5cm}x\in C\\
            x^+\in G(x), \hspace{.5cm}x\in D
            \end{array}\right.
            \end{equation}
        where the \textit{flow map}, $F: \mathbb{R}^n\rightrightarrows\mathbb{R}^n$ governs continuous flow of $x\in\mathbb{R}^n$, the \textit{flow set} $C\subset\mathbb{R}^n$ dictates where the continuous flow could occur. The \textit{jump map}, $G: \mathbb{R}^n\rightrightarrows\mathbb{R}^n$, governs discrete jumps of the state $x$, and the \textit{jump set} $D\subset\mathbb{R}^n$ defines where the discrete jumps are permitted. Note that the state $x\in\mathbb{R}^n$ could possibly include both continuous and discrete components. A \textit{hybrid arc} is a function $x: \textrm{dom} \:x\to\mathbb{R}^n$, where $\textrm{dom}\:x$ is a hybrid time domain and, for each fixed $j$, $t\mapsto x(t,j)$ is a locally absolutely continuous function on the interval $I_j=\{t: (t,j)\in\textrm{dom}\:x\}$.
%
\section{Synergistic families of potential functions on $SO(3)$}\label{synergistic}
  Given a finite index set $\mathcal{Q}\subset\mathbb{N}$, we let $\mathcal{C}^1\left(SO(3)\times\mathcal{Q},\mathbb{R}^+\right)$ denote the set of positive-valued continuously differentiable functions $\mathcal{U}: SO(3)\times\mathcal{Q}\to\mathbb{R}^+$, that is to say, for each $q\in\mathcal{Q}$, the map $R\mapsto\mathcal{U}(R,q)$ is continuous and differentiable on $SO(3)$.        %
            Additionally, for all $(R,q)\in SO(3)\times\mathcal{Q}$, $\nabla\mathcal{U}(R,q)\in T_RSO(3)$ denotes the gradient of $\mathcal{U}$, with respect to $R$, relative to the Riemannian metric \eqref{metric}. Let $\Psi_\mathcal{U}\subset SO(3)\times\mathcal{Q}$ denote the set of critical points of $\mathcal{U}$ where its gradient vanishes $\nabla\mathcal{U}(R,q)=0$. A function $\mathcal{U}\in\mathcal{C}^1\left(SO(3)\times\mathcal{Q},\mathbb{R}^+\right)$ is said to be  a \textit{potential function}
            with respect to the set $\mathcal{A}\subseteq SO(3)\times\mathcal{Q}$ if $\mathcal{U}(R,q)>0$ for all $(R,q)\notin\mathcal{A}$, and $\mathcal{U}(R,q)=0$, for all $(R,q)\in\mathcal{A}$. The set of all potential functions on $SO(3)\times\mathcal{Q}$ with respect to $\mathcal{A}$ is denoted as $\mathcal{P}(\mathcal{A})$, where a function $\mathcal{U}(R,q)\in\mathcal{P}(\mathcal{A})$ can be seen as a family of potential functions on $SO(3)$ encoded into a single function indexed by the variable $q$.
\begin{definition}\label{definition::synergism}
\cite{Mayhew2013}
  For a given finite index set $\mathcal{Q}\subset\mathbb{N}$, we let $\mathcal{A}=\{I\}\times\mathcal{Q}$ and $\mathcal{U}\in\mathcal{P}(\mathcal{A})$. The potential function $\mathcal{U}$ is said to be 
                \textit{centrally synergistic} if and only if there exist a constant $\delta>0$ such that
                \begin{align}\label{condition::synergistic}
                \bar{\delta}:=\underset{(R,q)\in\Psi_\mathcal{U}\setminus\mathcal{A}}{\mathrm{min}}&\left[\mathcal{U}(R,q)-\underset{p\in\mathcal{Q}}{\mathrm{min}}\;\mathcal{U}(R,p)\right]>\delta,
                                \end{align}
                where $\Psi_\mathcal{U}$ defines the set of all critical points of $\mathcal{U}$. The scalar $\bar\delta$ is referred to as the \textit{synergistic gap} of $\mathcal{U}$.
\end{definition}
The adjective \textit{``centrally''} refers to the fact that all the potential functions $R\mapsto\mathcal{U}(R,q)$ share the identity element $I$ as a critical point (non-centrally synergistic families were also proposed in \cite{Mayhew2013, lee2015tracking}). We may drop the adjective where not needed.  Condition \eqref{condition::synergistic} implies that at any given undesired critical point $(R,q)$ of $\mathcal{U}\in\mathcal{P}(\mathcal{A})$ there exists another point $(R,p)\in SO(3)\times\mathcal{Q}$ such that $\mathcal{U}(R,p)$ has a lower value than $\mathcal{U}(R,q)$. In the remainder of this paper, we let
            \begin{align*}
            \mathcal{A}:=\{I\}\times\mathcal{Q}.
            \end{align*}
 
It the recent literature, it was shown that once a synergistic family of potential functions on $SO(3)$ is obtained, a hybrid feedback controller that achieves global asymptotic stability results immediately follows \cite{mayhew2011synergistic, Mayhew2011, Mayhew2013, lee2015tracking}. The idea in \cite{mayhew2011synergistic} consists of stretching and compressing $SO(3)$ by applying the following transformation on $SO(3)\times\mathcal{Q}$
\begin{equation}\label{may_gam}
\Gamma(R,q)=\mathcal{R}_a(k_qP(R),u)R,
\end{equation}
where $u\in\mathbb{S}^2$ is a constant unit vector, $k_q\in\mathbb{R}$ is an indexed scalar gain, $P$ is a smooth positive definite function on $SO(3)$ with respect to $I$. By composing the map $\Gamma$ with an existing potential function, one can relocate the critical points while leaving the identity element a fixed point for all $q\in\mathcal{Q}$. Despite the originality of this approach, it was abandoned mainly due to the difficulty in finding an explicit expression of the synergistic gap.
\\
In this section, we build up from the ideas in \cite{mayhew2011synergistic} towards more generic constructions of central synergistic potential functions on $SO(3)$ via ``angular warping", while providing a thorough analysis of the synergism properties. Let $u\in\mathbb{S}^2$ be a fixed unit vector. Let us consider the map $\Gamma: SO(3)\times\mathcal{Q}\to SO(3)$ such that
\begin{equation}\label{Gamma_Ys}
\begin{array}{l}
\Gamma(R,q):=R\mathcal{R}_a(\theta_q(R),u),
\end{array}
\end{equation}
where $\theta_q:SO(3)\to\mathbb{R}$ is a real-valued function which is injective with respect to the index $q$. Note that the map in (\ref{may_gam}) uses a \textit{left} multiplication of $R$ by the additional rotation $\mathcal{R}_a(k_qP(R),u)$. In this work, however, we chose a \textit{right} multiplication of $R$ by $\mathcal{R}_a(\theta_q(R),u)$ where the map $\theta_q$ is to be designed later. This choice will allow us to express the control input directly in terms of vector measurements as done in Section \ref{Hybrid}.
\begin{lemma}\label{Lem_Gamma}
Let $u\in\mathbb{S}^2$ and $\mathcal{Q}\subset\mathbb{N}$. If the map $\theta_q(R): SO(3)\to\mathbb{R}$ is differentiable then the following properties hold:
\begin{enumerate}
\item The time derivative of $\Gamma(R,q)$ along the trajectories of $\dot R=R[\omega]_\times$ is given by
\begin{equation}\label{dGamma}
\frac{d}{dt}\Gamma(R,q)=\Gamma(R,q)[\Theta(R,q)\omega]_\times,
\end{equation}
where $
\Theta(R,q)=\mathcal{R}_a(\theta_q(R),u)^{\top}+2u\psi(R^{\top}\nabla\theta_q(R))^{\top}$.
\item If $\mathrm{det}(\Theta(R,q))\neq 0$ for all $(R,q)\in SO(3)\times\mathcal{Q}$, then the map $R\mapsto\Gamma(R,q)$ is everywhere a local \textit{diffeomorphism}. Moreover, if $V: SO(3)\to\mathbb{R}^+$ is a smooth positive definite function on $SO(3)$ with respect to $I$ and $\Gamma^{-1}(\{I\})=\{I\}\times\mathcal{Q}$ then $\mathcal{U}=V\circ\Gamma\in\mathcal{P}(\mathcal{A})$ with $\mathcal{A}=\{I\}\times\mathcal{Q}$ and the set of critical points of $\mathcal{U}$ is given by $\Psi_\mathcal{U}=\Gamma^{-1}(\Psi_V)$.
\end{enumerate}
\end{lemma}
Lemma \ref{Lem_Gamma} shows that, under some conditions on the transformation $\Gamma$, one can construct a new family of potential functions on $SO(3)$ by considering the composition of a \textit{basic} potential function on $SO(3)$ and the map $\Gamma$. In particular, it would be interesting to consider the modified trace function $V_A(R)=\mathrm{tr}(A(I-R))$ as the basic potential function due to its nice properties. The following technical lemma gives some of the useful properties of the potential function $V_A$.
\begin{lemma}\label{Lem_V_A}
Let $A=A^{\top}$ and $
V_A(R)=\mathrm{tr}(A(I-R)),
$
such that $W:=\mathrm{tr}(A)I-A$ is symmetric positive definite. Let $\{v_1,v_2,v_3\}$ be an orthonormal eigenbasis, where $v_i$ is a unit eigenvector associated to the eigenvalue $\lambda_i^A$. Then, for all $R\in SO(3)$, the following properties hold:
\begin{align}\label{gradV_A}
\nabla V_A(R)&=R\mathbb{P}_a(AR)\in T_RSO(3),\\\label{CritV_A}
\Psi_{V_A}&=\{I\}\cup\mathcal{R}_a(\pi,\mathcal{E}(A)),
\end{align}
where $\mathcal{E}(A)$ denotes the set of (real) unit eigenvectors of $A$. Moreover, for all $(\theta,u)\in\mathbb{R}\times\mathbb{S}^2$, one has
\begin{align}
\label{V_A_Q}
V_A(\mathcal{R}_a(\theta,u))&=2\sin^2(\theta/2)u^{\top}Wu,\\
V_A\left(\mathcal{R}_a(\pi,v)\mathcal{R}_a(\theta,u)\right)&=2\lambda^W-2\sin^2(\theta/2)\Delta(v,u),
\end{align}
where $\lambda^W$ denotes the eigenvalue of $W$ associated to the eigenvector $v\in\mathcal{E}(W)\equiv\mathcal{E}(A)$ and $\Delta(v,u)$ is computed as follows.
\begin{enumerate}
\item If $A=\lambda I_3$ $\left(\lambda_i^A=\lambda,\;\;i=1,2,3\right)$, then $\mathcal{E}(A)\equiv\mathbb{S}^2$ and
	\begin{equation*}
	\Delta(v,u)=\lambda^W\cos^2(\phi), \;\;\;\phi=\angle(u,v)
	\end{equation*}
\item If $A$ has two distinct eigenvalues $\lambda_1^A=\lambda_2^A\neq\lambda_3^A$, then $\mathcal{E}(A)=\left\{v_{12},v_3\right\}$, $v_{12}\in\mathrm{span}\{v_1,v_2\}\cap\mathbb{S}^2$, and
\begin{align}
\nonumber
\Delta(v_{12},u)&=(1-(u^\top v_3)^2)[\lambda_2^W-\lambda_3^W\sin^2(\phi)],\\\nonumber
\Delta(v_3,u)&=\left(\lambda_3^W-\lambda_2^W(1-(u^\top v_3)^2)\right),
\end{align}
such that $\phi=\angle(v_{12},u^\bot)$ and $u^\bot$ is the projection of $u$ on the plane $\mathrm{span}\{v_1,v_2\}$.
\item If $A$ has three distinct eigenvalues $0<\lambda_1^A<\lambda_2^A<\lambda_3^A$, then $\mathcal{E}(A)=\{v_1,v_2,v_3\}$ and
\begin{align*}
\Delta(v_1,u)&=\lambda_1^W-(u^\top v_2)^2\lambda_3^W-(u^\top v_3)^2\lambda_2^W,\\
\Delta(v_2,u)&=\lambda_2^W-(u^\top v_3)^2\lambda_1^W-(u^\top v_1)^2\lambda_3^W,\\
\Delta(v_3,u)&=\lambda_3^W-(u^\top v_1)^2\lambda_2^W-(u^\top v_2)^2\lambda_1^W.
\end{align*}
\end{enumerate}
\end{lemma}
In \cite{casau2014globally}, using the transformation $\Gamma$ defined in \eqref{may_gam}, the authors derived sufficient and necessary conditions such that the potential function $V_A\circ\Gamma$ is synergistic in the case where $A$ has distinct eigenvalues. In the following theorem we provide sufficient and necessary conditions, using our transformation $\Gamma$ defined in \eqref{Gamma_Ys}, such that the potential function $V_A\circ\Gamma$ is synergistic for an \textit{arbitrary} spectrum of $A$ and a \textit{general} angle function $\theta_q(\cdot)$.
\begin{theorem}\label{Th_synergism}
Let $A=A^\top$ such that $\mathrm{tr}(A)I-A$ is a symmetric positive definite matrix. Let $\mathcal{Q}\subset\mathbb{N}$ be an index set of \textit{finite} cardinality greater than or equal to $2$. Consider the transformation $\Gamma$ defined in \eqref{Gamma_Ys} and assume that $\mathrm{det}(\Theta(R,q))\neq 0$ for all $(R,q)\in SO(3)\times\mathcal{Q}$ and  $\Gamma^{-1}(\{I\})=\{I\}\times\mathcal{Q}$. The potential function
$$
\mathcal{U}(R,q):=V_A(\Gamma(R,q))=\mathrm{tr}\left(A(I-\Gamma(R,q))\right),
$$
is synergistic if and only if
\begin{equation}\label{Th_cond}
\Delta(v,u)>0,\;\;\;\;\;\textrm{for all}\;\;v\in\mathcal{E}(A),
\end{equation}
where $\Delta(v,u)$ is given in Lemma \ref{Lem_V_A}.
\end{theorem}
%
%
Theorem \ref{Th_synergism} provides necessary and sufficient conditions of synergism for the family of \textit{perturbed} modified trace functions $V_A(\Gamma(R,q))$. The condition $\Delta(v,u)>0,\;\forall v\in\mathcal{E}(A)$, imposes a constraint on the choice of the direction $u$ of the angular warping and the spectrum of the weighting matrix $A$. The following proposition discusses the feasibility of the synergy conditions of Theorem \ref{Th_synergism}. 
\begin{prop}\label{proposition::feasability}
Let $A=A^\top$ such that $\mathrm{tr}(A)I-A$ is a symmetric positive definite matrix. Let $\varrho_{ij}=(\lambda_i^A+(-1)^j\lambda_1^A)/(\lambda_3^A+\lambda_2^A)$. Then, the synergy condition (\ref{Th_cond}) is
\begin{itemize}
\item not satisfied if $A=\lambda I_3$ $\left(\lambda_i^A=\lambda,\;\;i=1,2,3\right)$ for all $u\in\mathbb{S}^2$.
\item satisfied if $A$ has two identical eigenvalues $0<\lambda_1^A=\lambda_2^A<\lambda_3^A$ and
$$
0<\varrho_{33}<(u^\top v_3)^2<1.
$$
\item satisfied if $A$ has distinct eigenvalues $0<\lambda_1^A<\lambda_2^A<\lambda_3^A$ and 
\begin{multline*}
-\varrho_{33}(u^\top v_1)^2+\varrho_{23}<(u^\top v_2)^2<-\varrho_{32}(u^\top v_1)^2+\varrho_{22},
\end{multline*}
\end{itemize}
where $v_i$ is the unit eigenvector of $A$ corresponding to the eigenvalue $\lambda_i^A$.
\end{prop}
Proposition \ref{proposition::feasability} suggests that a \textit{necessary} condition for synergism of $V_A(\Gamma(R,q))$ is that the weighting matrix $A$ must have at least two distinct eigenvalues. Furthermore, the direction $u$ of the angular warping should be carefully chosen with respect to the eigenvectors of $A$ in order for the potential function $V_A(\Gamma(R,q))$ to be synergistic. Later in Proposition \ref{proposition::maximization}, we will provide an \textit{optimal} choice of the unit vector $u\in\mathbb{S}^2$ that satisfies the feasibility conditions of Proposition \ref{proposition::feasability} while maximizing the synergistic gap.

When the synergism condition \eqref{Th_cond} is verified, it is important to explicitly compute the value of the synergistic gap required for the implementation of a synergistic hybrid controller. To do so, one needs to calculate the undesired critical points of $V_A(\Gamma(R,q))$ to evaluate the expression of the synergistic gap in \eqref{condition::synergistic}. These undesired critical points are obtained by solving equation \eqref{critical_point} for the unknown $R\in SO(3)$. Equation \eqref{critical_point}, along with the result of Lemma \ref{Lem_V_A}, yields
\begin{equation}\label{algebraic_VA}
V_A(R)=2\lambda^W-2\sin^2\left(\theta_{q}(R)/2\right)\Delta(v,u), \quad (R,q)\in\Psi_\mathcal{U}\setminus\mathcal{A}.
\end{equation}
In \cite{mayhew2011synergistic}, the warping angle was chosen as $\theta_q(R)=k_qV_A(R)$, thus the above equation leads to the following nonlinear algebraic equation for the unknown $V_A(R)$,
$$
V_A(R)=2\lambda^W-2\sin^2\left(k_qV_A(R)/2\right)\Delta(v,u).
$$
The above equation can be solved numerically but an explicit solution is hard to obtain. Instead, we propose the following choice of the warping angle $\theta_q: SO(3)\to (-\pi,\pi)$,
\begin{equation}\label{theta}
\theta_q(R)=2\arcsin\left(k_qV_A(R)\right),\;\;\;k_q\neq 0,
\end{equation}
that leads, in view of \eqref{algebraic_VA}, to a quadratic equation in $V_A(R)$ which can be solved to obtain
\begin{align}\label{V_crit}
V_A(R)=\frac{-1+\sqrt{1+16\lambda^Wk_q^2\Delta(u,v)}}{4k_q^2\Delta(u,v)}, \quad (R,q)\in\Psi_\mathcal{U}\setminus\mathcal{A}.
\end{align}
Once the value of $V_A(R)$ is obtained at the undesired critical points of $\mathcal{U}=V_A(\Gamma(R,q))$, one can compute these undesired rotations as follows
\begin{align}
\nonumber
R&=\mathcal{R}_a(\pi,v)\mathcal{R}_a(\theta_q(R),u)^{\top},\\\label{critical_point2}
\theta_q(R)&=2\arcsin\left(k_qV_A(R)\right), \quad (R,q)\in\Psi_\mathcal{U}\setminus\mathcal{A}.
\end{align}
Note that the scalar gain $k_q$ needs to be selected to ensure that $\theta_q(\cdot)$ is well defined for all $R\in SO(3)$. Also, one must make sure that the conditions of Lemma \ref{Lem_Gamma} are verified to guarantee that $V_A(\Gamma(R,q))$ is a suitable potential function on $SO(3)\times\mathcal{Q}$ with respect to $\{I\}\times\mathcal{Q}$. This is the purpose of the next proposition.
\begin{lemma}\label{Novelmap:prop}
Let $u\in\mathbb{S}^2$ be a unit vector and $A=A^\top$ such that $W=\mathrm{tr}(A)I-A$ is a symmetric positive definite matrix. Consider the transformation $\Gamma$ as defined in \eqref{Gamma_Ys} where $\theta_q(\cdot)$ is given by \eqref{theta}. If the scalar gain $k_q$ satisfies the inequality
\begin{equation}\label{k_max}
|k_q|<\bar{k}=\frac{1}{2\lambda_{\mathrm{max}}^W\sqrt{6-\max\{1,4\xi^2\}}},
\end{equation}
where $\xi=\lambda_{\min}^W/\lambda_{\max}^W$. Then $\Gamma^{-1}(\{I\})=\{I\}\times\mathcal{Q}$ and $\mathrm{det}\left(\Theta(R,q)\right)\neq 0$, for all $R\in SO(3)$ and $q\in\mathcal{Q}$. Furthermore, the gradient of the angle warping function $\theta_q$ is given by
\begin{equation}\label{dtheta}
\nabla\theta_q(R)=\frac{2k_qR\mathbb{P}_a(AR)}{\sqrt{1-k_q^2V_A^2(R)}}.
\end{equation}
\end{lemma}
As a consequence of Lemma \ref{Novelmap:prop}, if the scalars $k_q$ satisfies \eqref{k_max} then, by Lemma \ref{Lem_Gamma}, the composite function $\mathcal{U}=V_A\circ\Gamma$ is a suitable potential function on $SO(3)\times\mathcal{Q}$ with respect to $\mathcal{A}=\{I\}\times\mathcal{Q}$. Once the set of critical points for the potential function $\mathcal{U}=V_A(\Gamma(R,q))$ is determined from \eqref{critical_point2}, the synergistic gap defined in \eqref{condition::synergistic} can be evaluated. In the following theorem, we explicitly provide the expression of the synergistic gap of $\mathcal{U}$ in the case of $\mathcal{Q}=\{1,2\}$. Note that, one of the features of the angular warping approach is that a synergistic family can be generated using only two potential functions on $SO(3)$.

\begin{theorem}\label{Proposition::gap}
Let $u\in\mathbb{S}^2$ and $A=A^\top$ such that $W=\mathrm{tr}(A)I-A$ is a symmetric positive definite matrix. Let $\mathcal{Q}:=\{1,2\}$ and let $k_1=-k_2=k$, with $k$ satisfying condition (\ref{k_max}).  Consider the transformation $\Gamma$ as defined in \eqref{Gamma_Ys}, where $\theta_q(\cdot)$ is given by \eqref{theta}. Assume that $\Delta(u,v)>0$ is satisfied for all $v\in\mathcal{E}(A)$, then $\mathcal{U}=V_A(\Gamma(R,q))$ is synergistic with a gap $\bar\delta$ given by 
\begin{equation}\label{delta_exact}
\bar{\delta}=\underset{v\in\mathcal{E}(A)}{\mathrm{min}}\;\sigma(k,\lambda^W,\Delta(v,u))\geq \sigma(k,\underline{\lambda},\underline{\Delta})
\end{equation}
such that
\begin{align}\label{sigma}
\sigma(k,\lambda^W,\Delta(v,u))&=8k^2\bar V^2(1-k^2\bar V^2)\Delta(v,u)\\
\bar V&=\frac{-1+\sqrt{1+16\lambda^Wk^2\Delta(v,u)}}{4k^2\Delta(v,u)}
\end{align}
where $\underline\lambda=\underset{v\in\mathcal{E}(A)}{\mathrm{min}}\lambda^W$ and $\underline\Delta=\underset{v\in\mathcal{E}(A)}{\mathrm{min}}\Delta(v,u)$.
\end{theorem}
%
%
It is interesting to figure out the \textit{optimal} vector $u\in\mathbb{S}^2$ such that the synergistic gap given in (\ref{delta_exact}) is maximized while the condition of synergism $\underset{v\in\mathcal{E}(A)}{\mathrm{min}}\;\Delta(v,u)>0$ is verified. Since $\sigma(k,\lambda^W,\Delta(v,u))$ is a strictly increasing function of $\Delta(v,u)$ as shown in the proof of Theorem \ref{Proposition::gap}, we perform the following maximization with respect to $u\in\mathbb{S}^2$.
\begin{equation}\label{maximization}
\underset{
u\in\mathbb{S}^2}{\mathrm{max}}\;\underline\Delta=\underset{
u\in\mathbb{S}^2}{\mathrm{max}}\left(\underset{v\in\mathcal{E}(A)}{\mathrm{min}}\;\Delta(v,u)\right).
\end{equation}
\begin{prop}\label{proposition::maximization}
The unit vector $u\in\mathbb{S}^2$, solution of the maximization \eqref{maximization}, satisfies the following:
\begin{itemize}
\item if $A$ has two distinct eigenvalues $0<\lambda_1^A=\lambda_2^A<\lambda_3^A$,\\
$$
(u^\top v_3)^2=1-\frac{\lambda_2^A}{\lambda_3^A}.
$$
\item if $A$ has three distinct eigenvalues $0<\lambda_1^A<\lambda_2^A<\lambda_3^A$,\\ 
\begin{align*}
u^\top v_1=0,\;(u^\top v_2)^2=\frac{\lambda^A_2}{\lambda^A_2+\lambda^A_3},\;(u^\top v_3)^2=\frac{\lambda^A_3}{\lambda^A_2+\lambda^A_3},
\end{align*}
if $
\lambda^A_2\geq\lambda^A_1\lambda^A_3/(\lambda^A_3-\lambda^A_1)$. Otherwise, the optimal solution is
\begin{equation*}
\begin{split}
(u^\top v_i)^2=1-4\frac{\prod_{j\neq i}\lambda_j^A}{\sum_{j\neq k}\lambda_j^A\lambda_k^A},\;\;\;i\in\{1,2,3\}.
\end{split}
\end{equation*}
\end{itemize}
\end{prop}
Proposition \ref{proposition::maximization} gives an optimal choice of the angular warping direction $u\in\mathbb{S}^2$ that maximizes the synergistic gap. Moreover, it is straightforward to verify that this choice of $u$ satisfies the feasibility conditions of Proposition \ref{proposition::feasability}. Consequently, a complete construction of a synergistic potential function via angular warping with an explicit maximized synergistic gap has been achieved.
\section{Hybrid velocity-free attitude stabilization using vector measurements}\label{Hybrid}
We assume that the rigid body is equipped with sensors that provide measurements in the body-attached frame, denoted by $b_i\in\mathbb{R}^3$ of constant and known inertial vectors $r_i\in\mathbb{R}^3,\;i=1,2,...,n\geq 2$, satisfying the following assumption:
\begin{assum}\label{Assumption}
At least three vectors, among the $n$ inertial vectors, are not collinear.
\end{assum}
It should be noted that this Assumption~\ref{Assumption} is needed in our analysis and does not exclude the case where measurements of only two non-collinear inertial vectors are available, say $b_1$ and $b_2$ corresponding to the non-collinear inertial vectors $r_1$ and $r_2$. In this case, one can always construct a third vector $b_3=b_1\times b_2$ which corresponds to the measurement of $r_3 = r_1\times r_2$.
The rigid body rotational dynamics are governed by
\begin{align}
\label{kinematics}
\dot R&=R[\omega]_\times,\\\label{dynamics}
J\dot\omega&=[J\omega]_\times\omega+\tau,
\end{align}
where $R\in SO(3)$ represents the attitude, $\omega\in\mathbb{R}^3$ being the angular velocity of the rigid body expressed in the body-attached frame and $J\in\mathbb{R}^{3\times 3}$ is the constant inertia matrix of the rigid body. The control torque, expressed in the body frame, is denoted by $\tau\in\mathbb{R}^3$. \\
Our objective is to design a hybrid control input torque $\tau$, using only vector measurements, guaranteeing robust global asymptotic stabilization of the attitude $R\in SO(3)$ to a desired constant reference $R_d\in SO(3)$.\\
Let us define the following auxiliary dynamic system
\begin{equation}\label{aux}
\dot{\hat R}=\hat R[\beta]_\times,
\end{equation}
with an arbitrary initial condition $\hat R(0)\in SO(3)$ and a design variable $\beta\in\mathbb{R}^3$ to be defined later. Let $X_h=RY_h^{\top}$, $h\in\{1,2\}$, with $Y_1=\hat R$ and $Y_2=R_d$. The rotation matrix $X_1$ describes the discrepancy between the actual rigid body orientation and the orientation provided by the auxiliary system \eqref{aux}, and the rotation matrix $X_2$ describes the discrepancy between the actual rigid body orientation and the desired orientation. Let $\mathcal{Q}\subset\mathbb{N}$ be an index set of \textit{finite} cardinality and let our state variables be
$
X=(X_1,X_2)\in\mathcal{D}_X:=SO(3)\times SO(3),$ and $q=(q_1,q_2)\in\mathcal{D}_q:=\mathcal{Q}\times\mathcal{Q}$.
For $h\in\{1,2\}$, let us define $A_h=\sum_{i=1}^n\rho_{ih}r_ir_i^{\top}$, where $\rho_{ih}>0$ are some positive scalars . Assumption 1 ensures that $A_h,\;h=1,2$, are positive definite matrices. For some arbitrary unit vectors $u_h\in\mathbb{S}^2,\;h=1,2,$ and a set of arbitrary scalars $k_{hq}$ satisfying the conditions
\begin{align*}
|k_{hq}|<\frac{1}{2\lambda_{\mathrm{max}}^{W_h}\sqrt{6-\max\{1,4\xi_h^2\}}},\;\;\;\;h=1,2,\;q\in\mathcal{Q},
\end{align*}
with $\xi_h=\lambda_{\min}^{W_h}/\lambda_{\max}^{W_h}$ and $W_h:=\mathrm{tr}(A_h)I-A_h,\;h=1,2$.
We define the following maps $\Gamma_h:SO(3)\times\mathcal{Q}\to SO(3)$ such that
\begin{align*}
\Gamma_h(R,q)&=R\mathcal{R}_a(\theta_{hq}(R),u_h),\\
\theta_{hq}(R)&=2\arcsin(k_{hq}V_{A_h}(R)),\\
V_{A_h}(R)&=\mathrm{tr}\left(A_h(I-R)\right),
\end{align*}
for $h=1,2$. According to Section \ref{synergistic}, and in view of the above definition of the maps $\Gamma_h,\;h=1,2,$ one ensures that 
\begin{equation*}
\mathcal{U}_h(R,q):=V_{A_h}\circ\Gamma_h(R,q),
\end{equation*}
are both two potential functions on $SO(3)\times\mathcal{Q}$ with respect to $\{I\}\times\mathcal{Q}$.
%
We propose the following hybrid switching law for the control input $\tau$ and the input $\beta$ of the auxiliary system (\ref{aux})
\begin{equation}\label{hcontroller}
\begin{array}{ll}
\underset{(X,q)\in C}{\underbrace{\begin{array}{l}
\tau=-2\sum_{h=1}^2Y_h^\top\Theta_h(X_h,q_h)^\top\psi(A_h\Gamma_h(X_h,q_h)),\\
\beta=Y_1^\top\Theta_1(X_1,q_1)^\top\psi(A_1\Gamma_1(X_1,q_1))\\
\dot q=0
\end{array}}}
\\
\underset{(X,q)\in D}{\underbrace{\begin{array}{l}
\tau^+=\tau\\
\beta^+=\beta\\
q^+=g(X)
\end{array}}}
\end{array}
\end{equation}
where $\Theta_h(R,q)=\mathcal{R}_a(\theta_{hq}(R),u_h)^\top+\frac{4k_{hq}u_h\psi(A_hR)^\top}{\sqrt{1-k^2_{hq}V^2_{A_h}(R)}}$ and
\begin{multline*}
g(X)=\left\{(q_1,q_2)\in\mathcal{Q}\times\mathcal{Q}: q_h=\mathrm{arg}\underset{p\in\mathcal{Q}}{\mathrm{min}}\;\mathcal{U}_h(X_h,p),\right.\\\left.\;h=1,2\right\},
\end{multline*}
and the sets $C,D\subset\mathcal{D}_X\times\mathcal{D}_q$ are given by
\begin{multline*}
C:=\left\{(X,q)\in\mathcal{D}_X\times\mathcal{D}_q:\mu_1(X_1,q_1)\leq\delta_1 \;\mathrm{and}\right.\\\left.\;\mu_2(X_2,q_2)\leq\delta_2\right\}
\end{multline*}
\begin{multline*}
D:=\left\{(X,q)\in\mathcal{D}_X\times\mathcal{D}_q: \mu_1(X_1,q_1)\geq\delta_1 \;\mathrm{or}\right.\\\left.\;\mu_2(X_2,q_2)\geq\delta_2\right\}
\end{multline*}
such that $
\mu_h(X_h,q_h)=\mathcal{U}_h(X_h,q_h)-\underset{p\in\mathcal{Q}}{\mathrm{min}}\;\mathcal{U}_h(X_h,p), \;\;\;h=1,2.
$
The hybrid controller (\ref{hcontroller}) results in the closed-loop system
{\small
\begin{equation}\label{closed-loop}
\begin{array}{l}
\underset{(X,q)\in C
}{\underbrace{\begin{array}{l}
\dot X_1=X_1\left[Y_1\omega-\Theta_1(X_1,q_1)^\top\psi(A_1\Gamma_1(X_1,q_1))\right]_\times:=X_1[\omega_1]_\times,\\
\dot X_2=X_2[Y_2\omega]_\times:=X_2[\omega_2]_\times,\\
J\dot\omega=[J\omega]_\times\omega-2\sum_{h=1}^2Y_h^\top\Theta_h(X_h,q_h)^\top\psi(A_h\Gamma_h(X_h,q_h)),\\
\dot q=0
\end{array}}}\\
\underset{(X,q)\in D}{\underbrace{\begin{array}{l}
X^+=X\\
\omega^+=\omega\\
q^+=g(X)
\end{array}}}
\end{array}
\end{equation}}
Since $Y_1=X_1^{\top}X_2Y_2$ and $Y_2$ is constant, it is clear that the closed loop dynamics (\ref{closed-loop}) are autonomous. The goal of this hybrid controller is to ensure global asymptotic stability of the set
\begin{equation*}
\bar{\mathcal{A}}=\{(X,\omega,q)\in\mathcal{D}_X\times\mathbb{R}^3\times\mathcal{D}_q:\;X_1=X_2=I,\;\omega=0\}.
\end{equation*}

\begin{theorem}\label{theorem::hybrid}
Consider system \eqref{kinematics}-\eqref{dynamics} and the auxiliary system \eqref{aux} under the hybrid control law given in \eqref{hcontroller}. Assume that $n$ vector measurements $b_i$, corresponding to the inertial vectors $r_i$, $i=1,\cdots,n\geq 2$ are available, and Assumption \ref{Assumption} holds. If the potential function $\mathcal{U}_1$, respectively $\mathcal{U}_2$, is synergistic with gap exceeding $\delta_1$, respectively $\delta_2$, then the set $\bar{\mathcal{A}}$ is globally asymptotically stable for the closed-loop system (\ref{closed-loop}).
\end{theorem}
In practice, it is useful to explicitly express the control inputs in terms of the available vector measurements. Such measurements can be obtained, for instance, from an Inertial Measurement Unit (IMU) that typically includes an accelerometer and a magnetometer providing,
respectively, measurements of the gravitational field and Earth's magnetic field expressed in the body frame. The following Lemma shows that the terms involved in our hybrid control scheme can be directly expressed in terms of the available inertial vector measurement. 
 
\begin{prop}\label{Lem_mesurable}
The following relations hold:
\begin{align}\label{V_A0}
V_{A_h}(X_h)&=\frac{1}{2}\sum_{i=1}^{n}\rho_{ih}\left|\left|b_i-Y_h^\top r_i\right|\right|^2,\\
\label{psi0}
\psi\left(A_hX_h\right)&=\frac{1}{2}Y_h\sum_{i=1}^{n}\rho_{ih} (b_i\times Y_h^\top r_i),
\end{align}
for $h=1,2$. Furthermore, for $q\in\mathcal{Q}$ let
$$
\theta_{hq}(X_h)=2\arcsin(k_{hq}V_{A_h}(X_h)),
$$ and 
$\hat b_{ih}(q)=Y_h^{\top}\mathcal{R}_a(\theta_{hq}(X_h),u_h) r_i$, for $h=1,2,$
then
\begin{align}\label{V_Aq}
\mathcal{U}_h(X_h,q)&=\frac{1}{2}\sum_{i=1}^{n}\rho_{ih}\left|\left|b_i-\hat b_{ih}(q)\right|\right|^2,\\
\nonumber
\psi\left(A_h\Gamma_h(X_h,q)\right)&=\frac{1}{2}\mathcal{R}_a(\theta_{hq}(X_h),u_h)^{\top}Y_h\\
\label{psiq}
&\hspace{3cm}\sum_{i=1}^{n}\rho_{ih} (b_i\times\hat b_{ih}(q)).
\end{align}
\end{prop}
We devote the next section to prove all the Lemmas, Propositions and Theorems presented throughout the paper.
\section{Proofs}
\subsection{Proof of Lemma \ref{Lem_Gamma}}
First, using the product rule and the fact that $R[v]_\times R^{\top}=[Rv]_\times$ for all $R\in SO(3)$ and $v\in\mathbb{R}^3$, it is straightforward to show that if $\dot R_s=R_s[\omega_s]_\times$, $s=1,2$, then
\begin{equation}\label{dYs}
\frac{d}{dt}(R_1R_2)=R_1R_2\left[R_2^{\top}\omega_1+\omega_2\right]_\times.
\end{equation}
On the other hand, using \eqref{metric} and \eqref{id7}, one has
{\small $$\dot\theta_q(R)=d\theta_q(R)\cdot\dot R=\langle\nabla\theta_q(R),R[\omega]_\times\rangle_R=2\;\psi(R^{\top}\nabla\theta_q(R))^{\top}\omega,$$} which, using the fact that $\frac{d}{dt}e^{\theta(t)A}=e^{\theta(t)A}\dot\theta(t)A$, yields
\begin{equation}\label{dRt}
\begin{split}
\frac{d}{dt}\mathcal{R}_a(\theta_q(R),u)&=\mathcal{R}_a(\theta_q(R),u)\dot\theta_q(R)[u]_\times\\
													  &=\mathcal{R}_a(\theta_q(R),u)\left[2u\;\psi(R^{\top}\nabla\theta_q(R))^{\top}\omega\right]_\times.
\end{split}
\end{equation}
Since $\Gamma(R,q)=R\mathcal{R}_a(\theta_q(R),u)$ then, in view of (\ref{dYs}) and (\ref{dRt}), one obtains
$
\dot\Gamma(R,q)=\Gamma(R,q)\left[\Theta(R,q)\omega\right]_\times
$
where $\Theta(R,q)=\mathcal{R}_a(\theta_q(R),u)^{\top}+2u\;\mathrm{vex}(R^{\top}\nabla\theta_q(R))^{\top}$. 

For some $q\in\mathcal{Q}$, let us define the mapping $\mathcal{T}_q(R)=\Gamma(R,q)$. The time-derivative of the map $\mathcal{T}_q$ is nothing but the differential of $\mathcal{T}_q$ in the tangent direction $\xi=\dot R=R[\omega]_\times$. Replacing $\omega=\mathrm{vex}(R^{\top}\xi)$ in equation (\ref{dGamma}) shows that
\begin{equation*}
d\mathcal{T}_q(R)\cdot\xi=\mathcal{T}_q(R)[\Theta(R,q)\mathrm{vex}(R^{\top}\xi)]_\times,\;\;\;\xi\in T_RSO(3).
\end{equation*}
It is clear that when the inverse of the matrix $\Theta(R,q)$ exists for all $R\in SO(3)$ the map $d\\mathcal{T}_q(R)\cdot\xi$ is an isomorphism (bijective). In fact, for all $y:=d\mathcal{T}_q(R)\cdot\xi\in T_{\mathcal{T}_q(R)}SO(3)$, the inverse is explicitly given by
$$
\xi=R\left[\left(\Theta(R,q)\right)^{-1}\mathrm{vex}\left(\mathcal{T}_q(R)^{\top}y\right)\right]_\times\in T_RSO(3).
$$
Consequently, the inverse function theorem \cite{Moore1994} guarantees that $\mathcal{T}_q(R)$ is a local diffeomorphism for all $R\in SO(3)$. Furthermore, when $\Gamma$ is everywhere a local diffeomorphism on $SO(3)$ and $V: SO(3)\times\mathbb{R}^+$ is a smooth positive definite function, the composite function $\mathcal{U}=V\circ\Gamma$ is a \textit{differentiable positive} function on $SO(3)$. If in addition $\Gamma^{-1}(\{I\})=\mathcal{A}$ then $\mathcal{U}(R,q)=0$ if and only if $(R,q)\in\mathcal{A}$ which shows that $\mathcal{U}\in\mathcal{P}(\mathcal{A})$. In this case, the time-derivative of $\mathcal{U}=V\circ\Gamma$ can be computed as follows
{\small
\begin{equation}\label{dU}
\begin{split}
\dot{\mathcal{U}}(R,q)&=dV(\Gamma(R,q))\cdot\dot\Gamma(R,q)\\
				&=\langle\nabla V(\Gamma(R,q)),\dot\Gamma(R,q)\rangle_{\Gamma(R,q)}\\
			 &=2\;\psi\left(\Gamma(R,q)^{\top}\nabla V(\Gamma(R,q))\right)^{\top}\psi\left(\Gamma(R,q)^{\top}\dot\Gamma(R,q)\right)\\
			 &=2\;\psi\left(\Gamma(R,q)^{\top}\nabla V(\Gamma(R,q))\right)^{\top}\Theta(R,q)\omega,
			 \end{split}
\end{equation}}
where \eqref{metric} and (\ref{dGamma}) have been used. Consequently, since the matrix $\Theta(R,q)$ is full rank, the set of critical points of $\mathcal{U}$ is
{\small 
\begin{multline*}
\Psi_\mathcal{U}=\left\{(R,q)\in SO(3)\times\mathcal{Q}:\psi\left(\Gamma(R,q)^{\top}\nabla V(\Gamma(R,q))\right)=0\right\}.
\end{multline*}}
Also, the set of critical points of $V(R)$ is given by
$
\Psi_{V}:=\{R\in SO(3): \psi(R^{\top}\nabla V(R))=0\}.
$
This shows that $R\in\Psi_\mathcal{U}$ if and only if $\Gamma(R,q)\in\Psi_{V}$, which proves that
$\Psi_\mathcal{U}=\Gamma^{-1}(\Psi_V).$
\subsection{Proof of Lemma \ref{Lem_V_A}}
Recall from Section \ref{notation} that the gradient $\nabla V_A(R)$ can be computed using the differential of $V_A(R)$ in an arbitrary tangential direction $\xi=R\Omega\in T_RSO(3)$, such that
\begin{equation}\label{def_grad}
\begin{split}
DV_A(R)\cdot\xi=\langle\nabla V_A(R),R\Omega\rangle_R=\lls R^{\top}\nabla V_A(R),\Omega\ggs.
\end{split}
\end{equation}
On the other hand, one can compute $DV_A(R)\cdot\xi$ by direct differentiation pf $V_A(R)=\mathrm{tr}(A(I-R))$ as follows
\begin{equation*}
DV_A(R)\cdot\xi=-\mathrm{tr}(A\xi)=-\mathrm{tr}(AR\Omega)=\lls AR,\Omega\ggs.
\end{equation*}
Thus, one has $DV_A(R)\cdot\xi=\lls \mathbb{P}_a(AR),\Omega\ggs$ since $\Omega$ is skew symmetric. Therefore, in view of (\ref{def_grad}), yields (\ref{gradV_A}). The proof of (\ref{CritV_A}) is given in [\cite{mayhew2011synergistic}, Lemma $2$]. Now, let $Q=(\eta,\epsilon)\in\mathbb{Q}$ and $(\theta,u)\in\mathbb{R}\times\mathbb{S}^2$ be the quaternion and the angle-axis representation of the attitude matrix $R$, respectively. Using the Rodrigues formula (\ref{rod}) one obtains
\begin{equation*}
V_A(R)=\mathrm{tr}(A(-2[\epsilon]_\times^2-2\eta[\epsilon]_\times))=-2\mathrm{tr}(A[\epsilon]_\times^2),
\end{equation*}
where we used $\mathrm{tr}(A[\epsilon]_\times)=\lls A,[\epsilon]_\times\ggs=0$ since $A$ is symmetric. Also, using  $[\epsilon]_\times^2=-\|\epsilon\|^2I+\epsilon\epsilon^{\top}$ and $\epsilon^{\top}A\epsilon=\textrm{tr}(A\epsilon\epsilon^{\top})$, one obtains
\begin{equation*}
\begin{split}
V_A(R)&=2\mathrm{tr}(\epsilon^{\top}\epsilon A-A\epsilon\epsilon^{\top})=2\epsilon^{\top}W\epsilon,
								  \end{split}
\end{equation*}
which yields \eqref{V_A_Q} by noting that $\epsilon=\sin(\theta/2)u$. 
%
%
%

Let $Q_v=(\eta_v,\epsilon_{v})$ be the quaternion representation of the attitude $P_v:=\mathcal{R}_a(\pi,v)\mathcal{R}_a(\theta,u)$ such that $P_v=\mathcal{R}_Q(Q_v)$, for some $\theta\in\mathbb{R}$ and $u,v\in\mathbb{S}^2$. Making use of the quaternion multiplication rule (\ref{Qmultiply}) and (\ref{R1R2}), and the fact that $\mathcal{R}_a(\pi,v)$ and $\mathcal{R}_a\left(\theta,u\right)$ correspond to the quaternion $(0,v)$ and $(\cos(\theta/2),\sin(\theta/2)u)$, respectively, the quaternion vector part of $P_v$ is given by
\begin{equation}\label{e_v}
\epsilon_v=\cos(\theta/2)v+\sin(\theta/2)(v\times u).
\end{equation}
We let $\alpha_i:=u^\top v_i$ denotes the coordinates of $u$ in the eigenbasis $\{v_i\}_{i\in\{1,2,3\}}$.\\
1) $A=\lambda I_3$\\
In this case, using the fact that $v^{\top}(v\times u)=0$ and $\|v\times u\|^2=\sin^2(\phi)=1-\cos^2(\phi)$ where $\phi=\angle(u,v)$, one can compute the value of $V_A(P_v)$ from equation (\ref{V_A_Q}) and (\ref{e_v}) to obtain
\begin{equation*}
\begin{split}
V_A(P_v)&=2\lambda^W\epsilon_v^{\top}\epsilon_v\\
			  &=2\lambda^W\left[\cos^2(\theta/2)+\sin^2(\theta/2)\|v\times u\|^2\right]\\
			   &=2\lambda^W-2\sin^2(\theta/2)\Delta(v,u),
\end{split}
\end{equation*}
where $\Delta(v,u)=\lambda^W\cos^2(\phi)$.\\
2) $A$ has two distinct eigenvalues $\lambda_1^A=\lambda_2^A\neq\lambda_3^A$\\
Since $\{v_1v_2,v_3\}$ is an orthonormal basis, one can decompose $\mathbb{R}^3=E_3\oplus E_3^{\bot}$, where $E_3:=\mathrm{span}\{v_3\}$ and $E_3^{\bot}:=\mathrm{span}\{v_1,v_2\}$. Therefore, each vector $v\in\mathbb{R}^3$ is uniquely decomposed as $v=v^{\|}+v^{\bot}$ such that $v^{\|}\in E_3$ and $v^{\bot}\in E_3^{\bot}$. Consequently, $W(v_3\times u)=W(v_3\times u^{\bot})=\lambda_2^W(v_3\times u^{\bot})$ which, in view of (\ref{e_v}), yields
\begin{equation}\label{We1}
\begin{split}
W\epsilon_{v_3}&=W\left(\cos(\theta/2)v_3+\sin(\theta/2)(v_3\times u)\right)\\
						&=\lambda_3^W\cos(\theta/2)v_3+\lambda_2^W\sin(\theta/2)(v_3\times u^{\bot}).
\end{split}
\end{equation}
Now, making use of (\ref{We1}) and the fact that $v_3^{\top}(v_3\times u^{\bot})=0$ and $\|v_3\times u^{\bot}\|^2=\|u^{\bot}\|^2=\alpha_1^2+\alpha_2^2$, equation (\ref{V_A_Q}) yields
\begin{equation*}
\begin{split}
V_A(P_{v_3})&=2\epsilon_{v_3}^{\top}W\epsilon_{v_3}\\
					&=2\lambda_3^W\cos^2(\theta/2)+2\lambda_2^W\sin^2(\theta/2)(\alpha_1^2+\alpha_2^2)\\
					&=2\lambda_3^W-2\sin^2(\theta/2)\Delta(v_3,u).
\end{split}
\end{equation*}
Let $v_{12}\in E_3^{\bot}\cap\mathbb{S}^2,$ and $\phi=\angle{(v_{12}, u^{\bot})}$. In view of (\ref{e_v}), and using the fact that $u=u^{\|}+u^{\bot}$, the quaternion vector $\epsilon_{v_{12}}$ is decomposed as $\epsilon_{v_{12}}=\epsilon_{v_{12}}^{\|}+\epsilon_{v_{12}}^{\bot}$ where
\begin{align*}
\epsilon_{v_{12}}^{\|}&=\sin(\theta/2)(v_{12}\times u^{\bot})\in E_3,\\
\epsilon_{v_{12}}^\bot&=\cos(\theta/2)v_{12}+\sin(\theta/2)(v_{12}\times u^{\|})\in E_3^{\bot}.
\end{align*}
Making use of the identity $\|u\times v\|^2=\|u\|^2\|v\|^2\sin^2(\phi)$, such that $\phi=\angle(u,v)$, the norms of $\epsilon_{v_{12}}^{\|}$ and $\epsilon_{v_{12}}^\bot$ can be computed as follows
\begin{equation*}
\begin{split}
\|\epsilon_{v_{12}}^{\|}\|^2&=\sin^2(\theta/2)\|v_{12}\times u^{\bot}\|^2\\
									&=\sin^2(\theta/2)\|v_{12}\|^2\|u^{\bot}\|^2\sin^2(\phi)\\
					  &=\sin^2(\theta/2)(\alpha_1^2+\alpha_2^2)\sin^2(\phi),
					  \end{split}
\end{equation*}
and
\begin{equation*}
\begin{split}
\|\epsilon_{v_{12}}^\bot\|^2&=\cos^2(\theta/2)\|v_{12}\|^2+\sin^2(\theta/2)\|v_{12}\|^2\|u^{\|}\|^2\\
								   &=1-\sin^2(\theta/2)\left(1-\alpha_3^2\right).
\end{split}
\end{equation*}
where the fact that $\|v_{12}\|=1$ has been used. Therefore, since $W\epsilon_{v_{12}}=\lambda_3^W\epsilon_{v_{12}}^{\|}+\lambda_2^W\epsilon_{v_{12}}^\bot$, one obtains
\begin{equation*}
\begin{split}
V_A(P_{v_{12}})&=2\epsilon_{v_{12}}^{\top}W\epsilon_{v_{12}}\\
		&=2\lambda_3^W\|\epsilon_{v_{12}}^{\|}\|^2+2\lambda_2^W\|\epsilon_{v_{12}}^\bot\|^2\\
		&=2\lambda_2^W-2\sin^2\left(\theta/2\right)(\alpha_1^2+\alpha_2^2)[\lambda_2^W-\lambda_3^W\sin^2(\phi)]\\
		&=2\lambda_2^W-2\sin^2\left(\theta/2\right)\Delta(v_{12},u).
\end{split}
\end{equation*}
3) $A$ has three distinct eigenvalues $0<\lambda_1^A<\lambda_2^A<\lambda_3^A$\\
Using (\ref{a_cross_b}) one has $v_m\times u=\sum_{n,l}\varepsilon_{mnl}\alpha_nv_l$, which allows to write
\begin{equation*}
\begin{split}
W\epsilon_{v_m}&=W\left(\cos(\theta/2)v_m+\sin(\theta/2)(v_m\times u)\right)\\
				        &=\cos(\theta/2)\lambda_m^Wv_m+\sin(\theta/2)\sum_{n,l}\varepsilon_{mnl}\alpha_n\lambda_l^Wv_l,
\end{split}
\end{equation*}
where the fact that $Wv_m=\lambda_m^Wv_m$ for all $m\in\{1,2,3\}$ has been used to obtain the last equality. Now, in view of (\ref{V_A_Q}) and using the fact that $v_m^{\top}v_l=0$ for $m\neq l$, one gets
\begin{equation*}
\begin{split}
V_A(P_{v_m})&=2\epsilon_{v_m}^{\top}W\epsilon_{v_m}\\
					&=2\cos^2\left(\theta/2\right)\lambda_m^W+2\sin^2\left(\theta/2\right)\sum_{n,l}\varepsilon_{mnl}^2\alpha_n^2\lambda_l^W\\
					&=2\lambda_m^W-2\sin^2(\theta/2)\left[\lambda_m^W-\sum_{n,l}\varepsilon_{mnl}^2\alpha_n^2\lambda_l^W\right]\\
					&=2\lambda_m^W-2\sin^2(\theta/2)\Delta(v_m,u).
\end{split}
\end{equation*}	
\subsection{Proof of Theorem \ref{Th_synergism}}
According to Definition \ref{definition::synergism}, the potential function $\mathcal{U}$ is synergistic if 
\begin{equation*}
\mathcal{U}(R,q)-\underset{p\in\mathcal{Q}}{\mathrm{min}}\;\mathcal{U}(R,p)>0,\;\;\;\;\;\;\;\forall (R,q)\in\Psi_\mathcal{U}\setminus\mathcal{A}.
\end{equation*}
According to Lemma \ref{Lem_V_A}, the set of critical points for the function $V_A$ is $
\Psi_{V_A}=\{I\}\cup\mathcal{R}_a(\pi,v),\;v\in\mathcal{E}(A).
$
Invoking Lemma \ref{Lem_Gamma}, the set of critical points for $\mathcal{U}=V_A\circ\Gamma$ is given by
\begin{equation*}
\Psi_\mathcal{U}=\mathcal{A}\cup_{v\in\mathcal{E}(A)}\Gamma^{-1}(\mathcal{R}_a(\pi,v)).
\end{equation*}
Let $v\in\mathcal{E}(A)$ and let $(R,q)=\Gamma^{-1}(\mathcal{R}_a(\pi,v))$. Thus, $\Gamma(R,q)=\mathcal{R}_a(\pi,v)$, which in view of (\ref{Gamma_Ys}), yields
\begin{equation}\label{critical_point}
\begin{split}
R=\mathcal{R}_a(\pi,v)\mathcal{R}_a(\theta_q(R),u)^{\top}.
\end{split}
\end{equation}
Therefore, for a given $p\in\mathcal{Q}$, one has
\begin{equation*}
\begin{split}
\Gamma(R,p)&=R\mathcal{R}_a(\theta_p(R),u)\\
&=\mathcal{R}_a(\pi,v)\mathcal{R}_a(\theta_q(R),u)^{\top}\mathcal{R}_a(\theta_p(R),u)
\\
	&=\mathcal{R}_a(\pi,v)\mathcal{R}_a\left(\tilde\theta_{pq}(R),u\right),
\end{split}
\end{equation*}
where $\tilde\theta_{pq}(R)=\theta_p(R)-\theta_q(R)$ and identity (\ref{Rt1Rt2}) has been used. Consequently, invoking Lemma \ref{Lem_V_A}, one obtains
\begin{equation*}
\mathcal{U}(R,p)=V_A\left(\Gamma(R,p)\right)=2\lambda^W-2\sin^2\left(\tilde\theta_{pq}(R)/2\right)\Delta(v,u).
\end{equation*}
On the other hand, one has $\mathcal{U}(R,q)=V_A\circ\Gamma(R,q)=V_A(\mathcal{R}_a(\pi,v))=2v^{\top}Wv=2\lambda^W$. Therefore, at any undesired critical point $(R,q)\in\Psi_\mathcal{U}\setminus\mathcal{A}$, one has
\begin{equation}\label{mu_i}
\begin{split}
\mathcal{U}(R,q)-\underset{p\in\mathcal{Q}}{\mathrm{min}}\;\mathcal{U}(R,p)&=2\;\underset{p\in\mathcal{Q}}{\mathrm{max}}\;\sin^2\left(\tilde\theta_{pq}(R)/2\right)\Delta(v,u).
\end{split}
\end{equation}
Since $\theta_q(\cdot)$ is injective with respect to $q$, one has $p\neq q$ implies that $\tilde\theta_{pq}(R)\neq 0$. Consequently, in view of \eqref{mu_i}, the necessary and sufficient condition for $\mathcal{U}$ to be synergistic is that $\Delta(v,u)>0$ for all possible $v\in\mathcal{E}(A)$.
\subsection{Proof of Proposition \ref{proposition::feasability}}
Since $A$ is a symmetric matrix, it follows that the eigenvectors of $A$ can be chosen to form an orthonormal basis of $\mathbb{R}^{3\times 3}$, denoted $\{v_1,v_2,v_3\}$. We let $\alpha_i:=u^\top v_i$ define the coordinates of $u$ in the eigenbasis of $A$. We have three possible cases:
\begin{enumerate}
\item $A=\lambda I_3$ $\left(\lambda_i^A=\lambda,\;\;i=1,2,3\right)$\\
In this case, the set of unit eigenvectors of $A$ is $\mathcal{E}(A)\equiv\mathbb{S}^2$ and according to Lemma \ref{Lem_V_A}, one has $\Delta(v,u)=\lambda^W\cos^2(\phi),$ where $\phi=\angle(u,v)$ for all $v\in\mathcal{E}(A)$. However, for any chosen $u\in\mathbb{S}^2$ there exists $v\in\mathcal{E}(A)$ such that $\angle(u,v)=\pi/2$, or equivalently $\Delta(v,u)=0$. Therefore, the synergy condition of Theorem \ref{Th_synergism} can not be verified in this case.
\item $A$ has two distinct eigenvalues $\lambda_1^A=\lambda_2^A\neq\lambda_3^A$\\
In this case, the set of unit eigenvectors of $A$ is $\mathcal{E}(A)=\left\{v_{12},v_3\right\}$, where $v_{12}\in\mathrm{span}\{v_1,v_2\}\cap\mathbb{S}^2$. Thus, according to Lemma \ref{Lem_V_A}, the synergy condition of Theorem \ref{Th_synergism} is written as
\begin{equation*}
\begin{split}
\Delta(v_3,u)&=\left(\lambda_3^W-\lambda_2^W(\alpha_1^2+\alpha_2^2)\right)>0\\
\Delta(v_{12},u)&=(\alpha_1^2+\alpha_2^2)[\lambda_2^W-\lambda_3^W\sin^2(\phi)]>0,
\end{split}
\end{equation*}
for all $\phi=\angle(v_{12},u^\bot)$ such that $u^\bot$ represents the projection of $u$ on the plane $\mathrm{span}\{v_1,v_2\}$. Since $\sin^2(\phi)\leq 1$ and $\alpha_1^2+\alpha_2^2=1-\alpha_3^2$, the above conditions are equivalent to $0<\lambda_2^W(1-\alpha_3^2)<\lambda_3^W<\lambda_2^W$ which leads to
\begin{equation}\label{Feasibility:cond2}
0<\left(1-\frac{\lambda_3^W}{\lambda_2^W}\right)<\alpha_3^2<1.
\end{equation}
\item $A$ has three distinct eigenvalues $0<\lambda_1^A<\lambda_2^A<\lambda_3^A$\\
In this case, the set of unit eigenvectors of $A$ is $\mathcal{E}(A)=\{v_1,v_2,v_3\}$ and, according to Lemma \ref{Lem_V_A}, the synergy condition of Theorem \ref{Th_synergism} is written as
\begin{equation*}
\Delta(v_m,u)=\lambda_m^W-\sum_{n,l}\varepsilon_{mnl}^2\alpha_n^2\lambda_l^W>0,\;\;\;\forall m\in\{1,2,3\}.
\end{equation*}
Using the definition of the Levi-Cevita symbol given in section (\ref{Sec_attitude}), it is straightforward to verify that $\Delta(v_m,u)$ can explicitly be written as
\begin{equation}\label{Deltas}
\begin{array}{l}
\Delta(v_1,u)=\lambda_1^W-\alpha_2^2\lambda_3^W-\alpha_3^2\lambda_2^W,\\
\Delta(v_2,u)=\lambda_2^W-\alpha_1^2\lambda_3^W-\alpha_3^2\lambda_1^W,\\
\Delta(v_3,u)=\lambda_3^W-\alpha_2^2\lambda_1^W-\alpha_1^2\lambda_2^W.
\end{array}
\end{equation}
Since $u$ is a unit vector it verifies the unit constraint $\alpha_3^2=1-\alpha_1^2-\alpha_2^2$. Thus, equations (\ref{Deltas}) are rewritten as
\begin{align*}
\Delta(v_1,u)&=(\lambda_1^W-\lambda_2^W)+\alpha_2^2(\lambda_2^W-\lambda_3^W)+\alpha_1^2\lambda_2^W\\
\Delta(v_2,u)&=(\lambda_1^W-\lambda_3^W)\alpha_1^2+\lambda_1^W\alpha_2^2-(\lambda_1^W-\lambda_2^W)\\
\Delta(v_3,u)&=-\lambda_2^W\alpha_1^2-\lambda_1^W\alpha_2^2+\lambda_3^W.
\end{align*}

Note that, since the eigenvalues of the matrix $W=\mathrm{tr}(A)I-A$ are given by $\lambda_i^W=\mathrm{tr}(A)-\lambda_i^A$, it is obvious that $\lambda_1^W>\lambda_2^W>\lambda_3^W>0$ which implies that $\Delta(v_1,u)>0$. The conditions $\Delta(v_2,u),\Delta(v_3,u)>0$ are equivalent to $\underline{\chi}(\alpha_1^2)<\alpha_2^2<\bar{\chi}(\alpha_1^2)$ where
\begin{align*}
\underline{\chi}(\alpha_1^2)&=-\frac{(\lambda_1^W-\lambda_3^W)}{\lambda_1^W}\alpha_1^2+\frac{(\lambda_1^W-\lambda_2^W)}{\lambda_1^W}\\
\bar{\chi}(\alpha_1^2)&=-\frac{\lambda_2^W}{\lambda_1^W}\alpha_1^2+\frac{\lambda_3^W}{\lambda_1^W}.
\end{align*}
Now, since $\lambda_2^W+\lambda_3^W-\lambda_1^W=2\lambda_1^A>0$, it is easy to verify that $\lambda_3^W/\lambda_1^W>(\lambda_1^W-\lambda_2^W)/\lambda_1^W$ and $\lambda_3^W/\lambda_2^W>(\lambda_1^W-\lambda_2^W)/(\lambda_1^W-\lambda_3^W)$. A sketch of the feasible region is given in Figure \ref{feasible_reg} which shows that there always exist $\alpha_1^2$ and $\alpha_2^2$ such $\Delta(v_2,u), \Delta(v_3,u)>0$. Replacing $\lambda_i^W=\mathrm{tr}(A)-\lambda_i^A$ in the above inequality leads to the item 3 of Proposition \ref{proposition::feasability}. Finally, one can conclude that in the case where the matrix $A$ has \textit{arbitrary} distinct eigenvalues, one can always find a unit vector $u\in\mathbb{S}^2$ such that a synergistic family of potential functions (via angular warping) is constructed.
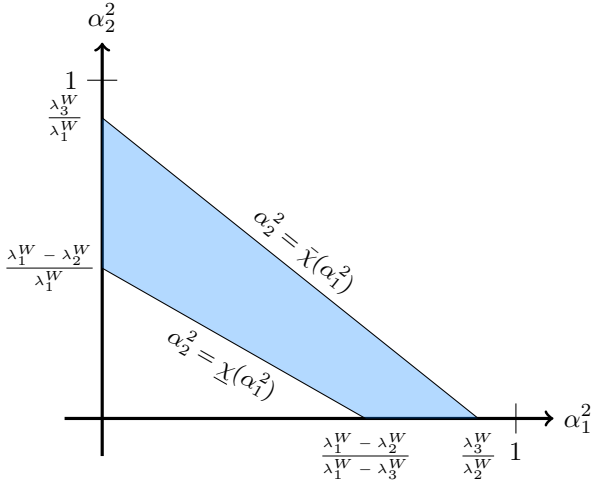
\begin{figure}[h!]
\centering
\begin{tikzpicture}

    \draw[very thick,->] (-.5,0) -- (6,0) node[right] {$\alpha_1^2$};
    \draw[very thick,->] (0,-.5) -- (0,5) node[above] {$\alpha_2^2$};

	\node at (-.5,4) {\tiny $\frac{\lambda_3^W}{\lambda_1^W}$};
	\node at (-.7,2) {\tiny $\frac{\lambda_1^W-\lambda_2^W}{\lambda_1^W}$};
		
	\node at (5,-0.5) {\tiny $\frac{\lambda_3^W}{\lambda_2^W}$};
	\node at (3.5,-0.5) {\tiny $\frac{\lambda_1^W-\lambda_2^W}{\lambda_1^W-\lambda_3^W}$};
	
    \fill[blue!50!cyan,opacity=0.3] (0,4) -- (0,2) --(3.5,0) -- (5,0) -- cycle;

    \draw (0,4) -- node[above,sloped] {\small $\alpha_2^2=\bar\chi(\alpha_1^2)$} (5,0);
    \draw (0,2) -- node[below,sloped] {\small $\alpha_2^2=\underline\chi(\alpha_1^2)$} (3.5,0);
	\draw (5.5,0.2) -- (5.5,-0.2)node[below] {$1$};
   \draw (0.2,4.5) -- (-0.2,4.5)node[left] {$1$};
\end{tikzpicture}
\caption{The feasible region of the synergy condition (\ref{Th_cond})} in the case where $A$ has distinct eigenvalues.
\label{feasible_reg}
\end{figure}
\end{enumerate}

\subsection{Proof of Lemma \ref{Novelmap:prop}}
In view of (\ref{V_A_Q}) one has $V_A(R)\leq 2\lambda_{\mathrm{max}}^W$ which, in view of (\ref{k_max}), implies that $|k_qV_A(R)|<1$. Therefore the function $\theta_q$ in (\ref{theta}) is well defined.
Let $(R,q)\in SO(3)\times\mathcal{Q}$ such that $\Gamma(R,q)=I$. This implies that $R=\mathcal{R}_a(\theta_q(R),u)^{\top}$ which, using (\ref{V_A_Q}) and the fact that $\sin(\theta_q(R)/2)=k_qV_A(R)$, yields
\begin{equation*}
\begin{split}
V_A(R)&=2\sin^2(\theta_q(R)/2)u^{\top}Wu=2k_q^2V_A^2(R)u^{\top}Wu,
\end{split}
\end{equation*}
or equivalently $(1-2k_q^2V_A(R)u^{\top}Wu)V_A(R)=0$. However, since $u\in\mathbb{S}^2$ and $|k_q|$ satisfies (\ref{k_max}), one has
$
2k_q^2V_A(R)u^{\top}Wu\leq 4k_q^2(\lambda_{\mathrm{max}}^W)^2<1/2.
$
Therefore, we must have $V_A(R)=0$ and thus $R=I$ since $V_A$ is positive definite on $SO(3)$. This shows that $\Gamma^{-1}(\{I\})=\{I\}\times\mathcal{Q}$.
Using the fact that the differentiability interval of $\arcsin$ is $(-1,1)$, and the fact that $|k_qV_A(R)|<1$, it is clear, by composition rule, that the function $\theta_q(R)=2\arcsin(k_qV_A(R))$ is differentiable on $SO(3)$. The gradient of $\theta_q(R)$ is computed using the chain rule and is given by
\begin{equation*}
\nabla\theta_q(R)=\frac{2k_q\nabla V_A(R)}{\sqrt{1-(k_qV_A(R))^2}}=\frac{2k_qR\mathbb{P}_a(AR)}{\sqrt{1-k_q^2V_A^2(R)}},
\end{equation*}
where (\ref{gradV_A}) has been used. Now, we prove the following inequality
\begin{equation}\label{ineq::psi}
\|\psi(AR)\|\leq\lambda_{\max}^W\min\{1,\sqrt{5-4\xi^2}/2\},
\end{equation}
where $\xi=\lambda_{\min}^W/\lambda_{\max}^W$. Let $(\eta,\epsilon)\in\mathbb{Q}$ be the quaternion representation of the attitude $R$. In view of \eqref{rod} and identity $[\epsilon]_\times^2=-\|\epsilon\|^2I+\epsilon\epsilon^{\top}$, one has
\begin{align*}
\mathbb{P}_a(AR)&=\frac{1}{2}(AR-R^\top A)\\
						   &=A\epsilon\epsilon^\top-\epsilon\epsilon^\top A+\eta A[\epsilon]_\times+\eta[\epsilon]_\times A\\
						   &=[\epsilon\times A\epsilon]_\times+\eta[W\epsilon]_\times,
\end{align*}
where equalities $yx^\top-xy^\top=[x\times y]_\times$ and $A^\top[x]_\times+A[x]_\times+[Ax]_\times=\mathrm{tr}(A)[x]_\times$, for all $x,y\in\mathbb{R}^3$, have been used. Consequently, one obtains 
$$
\psi(AR)=\epsilon\times A\epsilon+\eta W\epsilon=(\eta I-[\epsilon]_\times)W\epsilon.
$$
Therefore, 
\begin{equation}\label{inpsi1}
\begin{split}
\|\psi(AR)\|^2&=\epsilon^{\top}W(\eta I+[\epsilon]_\times)(\eta I-[\epsilon]_\times)W\epsilon\\
									   &=\epsilon^{\top}W(\eta^2I-[\epsilon]_\times^2)W\epsilon\\
									  &=\epsilon^{\top}W(I-\epsilon\epsilon^{\top})W\epsilon\\
									  &\leq (\lambda_{\max}^W)^2\|\epsilon\|^2-(\lambda_{\min}^W)^2\|\epsilon\|^4,
\end{split}
\end{equation}
which yields the fact that $\|\psi(AR)\|^2\leq (\lambda_{\max}^W)^2$ and also
\begin{align*}
\|\psi(AR)\|^2&	\leq (\lambda_{\max}^W)^2\|\epsilon\|^2\left(1-\|\epsilon\|^2+\|\epsilon\|^2-\xi^2\|\epsilon\|^2\right)\\
					&	\leq (\lambda_{\max}^W)^2\|\epsilon\|^2(1-\|\epsilon\|^2)+	(\lambda_{\max}^W)^2\|\epsilon\|^4(1-\xi^2)\\
					&\leq  (\lambda_{\max}^W)^2/4+	(\lambda_{\max}^W)^2(1-\xi^2),
\end{align*}
where we used the fact that $x^2(1-x^2)\leq 1/4$ for any $x\in[0 1]$. We just proved \eqref{ineq::psi}. Consequently,
\begin{equation}\label{nab_theta_ineq}
\begin{split}
\|\psi(R^{\top}\nabla\theta_q(R))\|&=\frac{2k_q\|\psi(AR)\|}{\sqrt{1-k_q^2V_A^2(R)}}\\
													&<\frac{2\lambda_{\max}^W\bar{k}\min\{1,\sqrt{5-4\xi^2}/2\}}{\sqrt{1-4\bar k^2(\lambda^W_{\mathrm{max}})^2}}\\
													&=\frac{\min\{1,\sqrt{5-4\xi^2}/2\}}{\sqrt{5-\max\{1,4\xi^2\}}}=\frac{1}{2},
	\end{split}
\end{equation}
where inequalities (\ref{ineq::psi}) and (\ref{k_max}) have been used along with $V_A(R)\leq 2\lambda_\mathrm{max}^W$. Now, making use of the matrix determinant lemma $\mathrm{det}(I+xy^{\top})=1+x^{\top}y$, for all $x,y\in\mathbb{R}^3$, one can compute the determinant of $\Theta(R,q)$ as follows:
\begin{align*}
\mathrm{det}\left(\Theta(R,q)\right)&=\mathrm{det}\left(\mathcal{R}_a(\theta_q(R),u)^{\top}\left[I+\right.\right.\\&\left.\left.\hspace{3.5cm}2u\psi(R^{\top}\nabla\theta_q(R))^{\top}\right]\right)\\
													 &=1+2u^{\top}\psi(R^{\top}\nabla\theta_q(R)),
\end{align*}
where $\mathcal{R}_a(\theta_q(R),u)^{\top}u=u$, $\mathrm{det}(AB)=\mathrm{det}(A)\mathrm{det}(B)$ and $\mathrm{det}(R)=1$ for all $R\in SO(3)$ have been used. Since $u$ is a unit vector and in view of (\ref{nab_theta_ineq}), it is obvious that the term $\left|2u^{\top}\psi(R^{\top}\nabla\theta_q(R))\right|$ is strictly less than $1$. Therefore $\mathrm{det}\left(\Theta(R,q)\right)\neq 0$ for all $R\in SO(3)$. 
\subsection{Proof of Theorem \ref{Proposition::gap}}
If $\Delta(u,v)>0$ is satisfied for all $v\in\mathcal{E}(A)$, then by Theorem \ref{Th_synergism}, the potential function $\mathcal{U}=V_A(\Gamma(R,q))$ is synergistic.  In order to compute the synergistic gap of $\mathcal{U}$, one needs to evaluate the expression in \eqref{mu_i}. Making use of the following trigonometric identity
\begin{equation*}
\arcsin(x)-\arcsin(y)=\arcsin\left(x\sqrt{1-y^2}-y\sqrt{1-x^2}\right),
\end{equation*}
it is clear that, for all $R\in SO(3)$, one has
\begin{equation*}
\begin{split}
\tilde\theta_{pq}(R)/2&=\theta_{p}(R)/2-\theta_{q}(R)/2\\
									&=\arcsin(k_pV_A(R))-\arcsin(k_qV_A(R))\\
									&=\arcsin\left(k_pV_A(R)\sqrt{1-k_q^2V_A^2(R)}-\right.\\
									&\left.\hspace{2cm}k_qV_A(R)\sqrt{1-k_p^2V_A^2(R)}\right).
\end{split}
\end{equation*}
Therefore, since $k_p=-k_q$ for $p\neq q$, one obtains
\begin{equation}\label{sin2}
\sin^2\left(\tilde\theta_{pq}(R)/2\right)=4k^2V_A^2(R)\left(1-k^2V_A^2(R)\right)^2.
\end{equation}
At any undesired critical point $(R,q)\in\Psi_{\mathcal{U}}\setminus\mathcal{A}$, in view of (\ref{mu_i}) and  (\ref{sin2}), one has
\begin{equation*}
\begin{split}
 \bar{\delta}:&=\underset{(R,q)\in\Psi_\mathcal{U}\setminus\mathcal{A}}{\mathrm{min}}\left[\mathcal{U}(R,q)-\underset{p\in\mathcal{Q}}{\mathrm{min}}\;\mathcal{U}(R,p)\right]\\
&=2\;\underset{(R,q)\in\Psi_\mathcal{U}\setminus\mathcal{A}}{\mathrm{min}}\left[\underset{p\in\mathcal{Q}}{\mathrm{max}}\;\sin^2\left(\tilde\theta_{pq}(R)/2\right)\Delta(v,u)\right]\\
 			 &=8k^2\bar V^2(1-k^2\bar V^2)\Delta(v,u),
\end{split}
\end{equation*}
where we used the fact that $V_A(R)=\bar V$ at $(R,q)\in\Psi_{\mathcal{U}}\setminus\mathcal{A}$ (see equation (\ref{V_crit})). By direct differentiation of \eqref{delta_exact} with respect to $\Delta(v,u)$, one obtains
\begin{multline*}
\frac{\partial\sigma}{\partial\Delta}=\frac{8k^2\bar V^2}{1+4k^2\Delta\bar V}\left\{1-k^2\bar V^2+\frac{\bar V}{2\Delta}\left(1+8\lambda^W k^2\Delta\right.\right.\\\left.\left.-\sqrt{1+16\lambda^Wk^2\Delta}\right)\right\}.
\end{multline*}
Now, using the fact that $1+x\geq\sqrt{1+2x}$ for all $x\geq 0$ and $|k\bar V|<1/\sqrt{2}$, one concludes that $\partial\sigma/\partial\Delta>0$. Also, one can verify that
\begin{equation*}
\frac{\partial\sigma}{\partial\lambda^W}=\frac{16k^2\Delta \bar V}{1+4k^2\Delta\bar V}(1-2k^2\bar V^2)>0.
\end{equation*}
Consequently, one has
$
\bar{\delta}\geq\sigma\left(k,\underline\lambda,\underline\Delta\right),
$
such that $\underline\lambda=\underset{v\in\mathcal{E}(A)}{\mathrm{min}}\lambda^W$ and $\underline\Delta=\underset{v\in\mathcal{E}(A)}{\mathrm{min}}\Delta(v,u)$.
\subsection{Proof of Proposition \ref{proposition::maximization}}
The term $\Delta(v,u)$ depends explicitly on the square of the coefficients $\alpha_i=u^\top v_i$, coordinates of $u$ in the eigenbasis $\{v_1,v_2,v_3\}$. Therefore, the maximization \eqref{maximization} can be performed with respect to $\alpha_i^2$ subject to the unit constraint $\alpha_1^2+\alpha_2^2+\alpha_3^2=1$. Depending on the spectrum of $A$ and excluding the case where $A=\lambda I$ (since no synergistic family can be constructed), we have two possible cases.
\begin{enumerate}
\item $A$ has two distinct eigenvalues $\lambda_1^A=\lambda_2^A\neq\lambda_3^A$\\
According to (\ref{Feasibility:cond2}), our maximization problem is transformed into finding $(1-\lambda_3^W/\lambda_2^W)<\alpha_3^2<1$ such that
\begin{multline*}
\underset{v\in\mathcal{E}(A)}{\mathrm{min}}\;\Delta(v,u)=\mathrm{min}\left\{\left(\lambda_3^W-\lambda_2^W(1-\alpha_3^2)\right),\right.\\\left.\left((\lambda_2^W-\lambda_3^W)(1-\alpha_3^2)\right)\right\},
\end{multline*}
is maximized.
\begin{figure}[h!]
\centering
\begin{tikzpicture}

    \draw[very thick,->] (-.5,0) -- (6,0) node[right] {$\alpha_3^2$};
    \draw[very thick,->] (0,-3) -- (0,3) node[above] {};

	\node at (-.7,2.5) {\tiny $\lambda_2^W-\lambda_3^W$};
	\node at (-.7,-2.5) {\tiny $\lambda_3^W-\lambda_2^W$};

	\node at (3,-0.5) {\tiny $1-\frac{\lambda_3^W}{\lambda_2^W}$};
	
    \fill[gray!50!white,opacity=0.3] (0,-3) -- (0,3) --(3.1,3) -- (3.1,-3) -- cycle;

    \draw (0,2.5) -- node[above,sloped] {\small $ $} (5.5,0);
    \draw (0,-2.5) -- node[below,sloped] {\small $ $} (5.5,2);
    \draw[very thick, green] (3.1,0)--(3.92,0.7) -- (5.5,0);
    \node at (5,0.7) {\tiny $\underset{v\in\mathcal{E}(A)}{\mathrm{min}}\;\Delta(v,u)$};

	\draw (5.5,0.2) -- (5.5,-0.2)node[below] {$1$};
\end{tikzpicture}
\caption{}
\label{Maximazation:plot}
\end{figure}
Figure (\ref{Maximazation:plot}) gives a plot of $\underset{v\in\mathcal{E}(A)}{\mathrm{min}}\;\Delta(v,u)$ with respect to $\alpha_3^2$ (green line). It is obvious that the maximum (top vertex of the green triangle) is attained at the intersection
$
\lambda_3^W-\lambda_2^W(1-\alpha_3^2)=(\lambda_2^W-\lambda_3^W)(1-\alpha_3^2),
$ which leads to
$$
\alpha_1^2=\frac{2(\lambda_2^W-\lambda_3^W)}{2\lambda_2^W-\lambda_3^W}=1-\frac{\lambda_2^A}{\lambda_3^A}.
$$
\item $A$ has three distinct eigenvalues $0<\lambda_1^A<\lambda_2^A<\lambda_3^A$\\
In this case, the maximization problem is formulated as
\begin{equation}\label{Maximization:case3}
\underset{
\alpha_i^2}{\mathrm{max}}\;\underline\Delta=\underset{
\alpha_i^2}{\mathrm{max}}\left(\underset{m\in\left\{1,2,3\right\}}{\mathrm{min}}\;\left(\lambda_m^W-\sum_{n,l}\varepsilon_{mnl}^2\alpha_n^2\lambda_l^W\right)\right).
\end{equation}
subject to
$
\sum_{i=1}^3\alpha_i^2=1$ and $
\underline\Delta>0$. The above optimization problem can be formulated as a standard linear programming problem as follows:
\begin{equation*}
\begin{array}{lrl}
\textrm{maximize}\quad&\;\hspace{-3cm}x_4&\\
\textrm{subject to}\quad&x_4+x_2\lambda_3^W+x_3\lambda_2^W-\lambda_1^W&\leq 0,\\
&x_4+x_3\lambda_1^W+x_1\lambda_3^W-\lambda_2^W&\leq 0,\\
&x_4+x_1\lambda_2^W+x_2\lambda_1^W-\lambda_3^W&\leq 0,\\
&x_1+x_2+x_3&=1,\\
&x_1,x_2,x_3&\geq 0,\\
&x_4&>0.
\end{array}
\end{equation*}
where the variables $x_j=\alpha_j^2$ for $j\in\{1,2,3\}$. Consequently, one can use the simplex algorithm \cite{robert1996linear} to solve this optimization problem, leading to the following result: If the condition \begin{equation}\label{cond_lambdas}
\lambda^A_2\geq\frac{\lambda^A_1\lambda^A_3}{\lambda^A_3-\lambda^A_1}
\end{equation}
is satisfied, then the optimal solution to the maximization in (\ref{Maximization:case3}) is given by
\begin{align*}
\alpha_1^2=0,\;\;\;\alpha_2^2=\frac{\lambda^A_2}{\lambda^A_2+\lambda^A_3},\;\;\;\alpha_3^2=\frac{\lambda^A_3}{\lambda^A_2+\lambda^A_3}
\end{align*}
and the maximum is $\underset{
\alpha_j^2}{\mathrm{max}}\;\underline\Delta=\lambda^A_1$. Otherwise, the optimal solution is
\begin{equation*}
\begin{split}
\alpha_i^2=1-4\frac{\prod_{j\neq i}\lambda_j^A}{\sum_{j\neq k}\lambda_j^A\lambda_k^A},\;\;\;i\in\{1,2,3\}
\end{split}
\end{equation*}
with the maximum being
$
\underset{
\alpha_j^2}{\mathrm{max}}\;\underline\Delta=4\frac{\prod_{j=1}^{3} \lambda_j^A}{\sum_{j\neq k}\lambda_j^A\lambda_k^A}
$
\end{enumerate}
\subsection{Proof of Theorem \ref{theorem::hybrid}}
Let us define the following sets:
\begin{align*}
\mathcal{A}_1&:=\{(X,q)\in\mathcal{D}_X\times\mathcal{D}_q:\;X_1=I\},\\
\mathcal{A}_2&:=\{(X,q)\in\mathcal{D}_X\times\mathcal{D}_q:\;X_2=I\},\\
\mathcal{X}_1&:=\{(X,q)\in\mathcal{D}_X\times\mathcal{D}_q:\;(X_1,q_1)\in\Psi_{\mathcal{U}_1}\},\\
\mathcal{X}_2&:=\{(X,q)\in\mathcal{D}_X\times\mathcal{D}_q:\;(X_2,q_2)\in\Psi_{\mathcal{U}_2}\},\\
C_1&:=\left\{(X,q)\in\mathcal{D}_X\times\mathcal{D}_q:\;\mu_1(X_1,q_1)\leq\delta_1 \right\},\\
C_2&:=\left\{(X,q)\in\mathcal{D}_X\times\mathcal{D}_q:\;\mu_2(X_2,q_2)\leq\delta_2 \right\}.
\end{align*}
Assume that $\mathcal{U}_h$ is synergistic with gap exceeding $\delta_h$, for $h=1,2$. Then, in view of (\ref{condition::synergistic}), one has
\begin{equation*}
0<\delta_h<\underset{(X_h,q_h)\in\Psi_{\mathcal{U}_h}\setminus\mathcal{A}}{\mathrm{min}}\mu_h(X_h,q_h).
\end{equation*}
Therefore for each pair $(X,q)\in\mathcal{X}_h\setminus\mathcal{A}_h$, one obtains
$
\mu_h(X_h,q_h)>\delta_h,
$
which implies, in view of the definition of the set $C_h$, that one has
\begin{equation}\label{CX}
C_h\cap\mathcal{X}_h=\mathcal{A}_h,
\end{equation}
where we used the fact that $\mathcal{A}_h$ is entirely contained in $C_h$.
Consider the Lyapunov function candidate
\begin{equation*}
\mathcal{V}(X,\omega,q)=\sum_{h=1}^2\mathcal{U}_h(X_h,q_h)+\frac{1}{2}\omega^{\top}J\omega.
\end{equation*}
Since $\mathcal{U}_h$ are two potential functions on $SO(3)\times\mathcal{Q}$ with respect to $\{I\}\times\mathcal{Q}$ and $J$ is positive definite, it follows that $\mathcal{V}$ is positive definite on $\mathcal{D}_X\times\mathbb{R}^3\times\mathcal{D}_q$ with respect to $\bar{\mathcal{A}}$.
In view of (\ref{dU}), the time derivative of $\mathcal{U}_h(X_h,q_h)$ along the trajectory $\dot X_h=X_h[\omega_h]_\times$ is given by
\begin{equation}\label{dUh}
\begin{split}
\dot{\mathcal{U}}_h(X_h,q_h)&=2\;\psi\left(\Gamma_h(X_h,q_h)^{\top}\nabla V_{A_h}(\Gamma_h(X_h,q_h))\right)^{\top}\\
										   &\hspace{5cm}\Theta_h(X_h,q_h)\omega_h\\
							&=2\psi\left(A_h\Gamma_h(X_h,q_h)\right)^{\top}\Theta_h(X_h,q_h)\omega_h,						
\end{split}
\end{equation}
where we used $\nabla V_A(R)=R\mathbb{P}_a(AR)$ for all $R\in SO(3)$. Therefore, making use of the above result and (\ref{closed-loop}), the change in $\mathcal{V}$ along the continuous flows of $C$ is given by
\begin{equation*}
\begin{split}
\dot{\mathcal{V}}(X,\omega,q)&=2\sum_{h=1}^2\omega_h^\top\Theta_h(X_h,q_h)^\top\psi(A_h\Gamma_h(X_h,q_h))+\\
&\hspace{-1cm}\omega^{\top}\left([J\omega]_\times\omega-2\sum_{h=1}^2Y_h^\top\Theta_h(X_h,q_h)^\top\psi(A_h\Gamma_h(X_h,q_h))\right)\\
&=-2\|\Theta_1(X_1,q_1)^\top\psi(A_1\Gamma_1(X_1,q_1))\|^2\leq 0.
\end{split}
\end{equation*}
Thus $\mathcal{V}$ is non-increasing along the flows of (\ref{closed-loop}). Moreover, for any $(X,q)\in D$ and $s\in g(X)$, one has
\begin{equation*}
\begin{split}
\mathcal{V}(X,\omega,q)-\mathcal{V}(X,\omega,s)&=\sum_{h=1}^2\left[\mathcal{U}_h(X_h,q_h)-\underset{p\in\mathcal{Q}}{\mathrm{min}}\;\mathcal{U}_h(X_h,p)\right]\\
&=\sum_{h=1}^2\mu_h(X_h,q_h)\\&\geq\mathrm{min}\{\delta_1,\delta_2\}>0,
\end{split}
\end{equation*}
which shows that $\mathcal{V}$ is strictly decreasing over the jumps of (\ref{closed-loop}). Using [\cite{Teel2007}, Theorem 7.6], it follows that $\bar{\mathcal{A}}$ is stable. Moreover, applying the invariance principle for hybrid systems given in [\cite{Teel2007}, Theorem 4.7], one can conclude that any solution must converge to the largest invariant set contained in
\begin{multline*}
\mathcal{I}=\left\{(X,\omega,q)\in\mathcal{D}_X\times\mathbb{R}^3\times\mathcal{D}_q: (X,q)\in C, \right.\\\left.\Theta_1(X_1,q_1)^\top\psi(A_1\Gamma_1(X_1,q_1))=0\right\}.
\end{multline*}
It follows, in view of (\ref{dUh}), that for all $(X,\omega,q)\in\mathcal{I}$, one has $(X_1,q_1)\in\Psi_{\mathcal{U}_1}$.
Consequently, the set $\mathcal{I}$ can be rewritten as
$$\mathcal{I}=\left\{(X,\omega,q)\in\mathcal{D}_X\times\mathbb{R}^3\times\mathcal{D}_q: (X,q)\in C\cap\mathcal{X}_1\right\}.$$
Moreover, from (\ref{CX}), one has $C_1\cap\mathcal{X}_1=\mathcal{A}_1$  and hence
$
C\cap\mathcal{X}_1=(C_1\cap C_2)\cap\mathcal{X}_1=C_2\cap\mathcal{A}_1,
$
where we used the fact that $C=C_1\cap C_2$. Since the solutions converge to $C_2\cap\mathcal{A}_1$, it is clear that $X_1\to I$ which leads to $\dot X_1\to 0$. Hence, one can conclude from (\ref{closed-loop}) that $\omega\to 0$. Since $\omega\equiv 0$, it follows from (\ref{dynamics}) that $\tau$ must converge to $0$. Using this last fact, together with the fact that $\beta\equiv 0$, one can conclude from (\ref{hcontroller}) that
$$
\Theta_2(X_2,q_2)^\top\psi(A_2\Gamma_2(X_2,q_2))=0.
$$
Again using (\ref{dUh}), one has $(X_2,q_2)\in\Psi_{\mathcal{U}_2}$. Therefore, the solutions must converge to
$
C_2\cap\mathcal{A}_1\cap\mathcal{X}_2=\mathcal{A}_1\cap\mathcal{A}_2=\bar{\mathcal{A}},
$
where the fact that $C_2\cap\mathcal{X}_2=\mathcal{A}_2$ has been used. Finally, the set $\bar{\mathcal{A}}$ is globally attractive and stable which shows that $\bar{\mathcal{A}}$ is globally asymptotically stable.
\subsection{Proof of Proposition \ref{Lem_mesurable}}
Let us define the attitude error $\tilde P:=RP^{\top}$, for some $P\in SO(3)$. For $h\in\{1,2\}$, making use of the identity $u^{\top}Au=\mathrm{tr}(uu^{\top}A)$, one obtains
\begin{equation*}
\begin{split}
\frac{1}{2}\sum_{i=1}^{n}\rho_{ih}\|b_i-P^{\top}r_i\|^2&=\frac{1}{2}\sum_{i=1}^{n}\rho_{ih}r_i^{\top}(I-\tilde P^{\top})(I-\tilde P)r_i\\
				   											 &=\sum_{i=1}^{n}\rho_{ih}\mathrm{tr}(r_ir_i^{\top}(I-\tilde P))=V_{A_h}(\tilde P).
\end{split}
\end{equation*}
Consequently, equation (\ref{V_A0}), respectively \eqref{V_Aq}, is obtained by substituting $P$ for $Y_h$, respectively $P$ for $\mathcal{R}^{\top}_a(\theta_{hq}(X_h),u_h)Y_h$. Furthermore, one has
\begin{equation}\label{P1}
\begin{split}
\mathbb{P}_a(A_h\tilde P)&=\frac{1}{2}\sum_{i=1}^{n}\rho_{ih}\left[r_ir_i^{\top}\tilde{P}-\tilde P^{\top}r_ir_i^{\top}\right]\\
&=\frac{1}{2}\sum_{i=1}^{n}\rho_{ih}P\left[P^\top r_ir_i^{\top}R-R^{\top}r_ir_i^{\top}P\right]P^\top\\
&=\frac{1}{2}\sum_{i=1}^{n}\rho_{ih}\left[P(b_i\times P^{\top}r_i)\right]_\times,
\end{split}
\end{equation}
where we used the following property (see \cite{Shuster1993}):
\begin{equation*}
R[yx^{\top}-xy^{\top}]R^{\top}=R[x\times y]_\times R^{\top}=[R(x\times y)]_\times.
\end{equation*}
Taking the $\mathrm{vex}$ operator on both sides of (\ref{P1}) and substituting $P$ for $Y_h$, respectively $P$ for $\mathcal{R}^{\top}_a(\theta_{hq}(X_h),u_h)Y_h$, yields equation (\ref{psi0}), respectively equation \eqref{psiq}.

\section{Simulation results}\label{Simulation}
In this section we illustrate the procedure to follow for the implementation of the hybrid scheme derived in Section \ref{Hybrid}. We, then, compare between the smooth feedback law, proposed in \cite{TAC2013}, and the hybrid feedback of Section \ref{Hybrid}. Let $\mathcal{Q}=\{1,2\}$ and let the scalars $k_{h1}=-k_{h2}=k_h$, for  $h=1,2$ where $k_h$ has a specific value to be tuned later. Let also $A_1=\mathrm{diag}([1,3,5])$ and $A_2=\mathrm{diag}([0.1,0.3,0.5])$. The value of $k_h$ is chosen to verify inequality (\ref{k_max}). We picked $k_1=0.03$ and $k_2=0.3$. Moreover, with this choice of matrices $A_1$ and $A_2$, one can see that condition (\ref{cond_lambdas}) is verified for both matrices. Therefore, the optimal choice of the rotation vectors $u_1$ and $u_2$ is given by
\begin{equation*}
u_1^{\top}=[0,\sqrt{3/8},\sqrt{5/8}],\;\;\;u_2^{\top}=[0,\sqrt{3/8},\sqrt{5/8}].
\end{equation*}
The maximum values of the hysteresis gaps $\bar\delta_1$ and $\bar\delta_2$ can be computed, using the result of Theorem \ref{Proposition::gap}. Therefore, it is sufficient to pick $\delta_1=0.5<\bar\delta_1$ and $\delta_2=0.05<\bar\delta_2$ to implement the switching conditions of the set $C$ and the set $D$. Once all the parameters have been designed, the hybrid controller (\ref{hcontroller}) can be implemented. We recall the following smooth feedback law proposed in the paper \cite{TAC2013}
\begin{equation}\label{smooth}
\begin{split}
\tau&=-\sum_{h=1}^2\sum_{i=1}^{n}\rho_{ih} (b_i\times Y_h^\top r_i),\\
\beta&= \sum_{i=1}^{n}\rho_{i1} (b_i\times Y_1^\top r_i).
\end{split}
\end{equation}
Both controllers, hybrid and smooth, were implemented in Simulink. The desired rotation as well as the initial condition $\hat R(0)$ for the auxiliary system were both chosen equal to the identity matrix, \textit{i.e.,} $R_d=\hat R(0)=I_{3\times 3}$. The inertia matrix has been taken as $J=\mathrm{diag}\left([1,1,2]\right)$ and the inertial vectors as $r_i=e_i$, for $i=1,2,3$. The performance of the controllers was evaluated by means of the \textit{normilized} error\footnote{This error can be shown to be equal to $e(X_h)=\|\epsilon_h\|_2$ where $\epsilon_h$ is the quaternion vector part corresponding to the orientation described by $X_h\in SO(3)$.}
$$
e(X_h):=\frac{1}{\sqrt{8}}\|I-X_h\|_F.
$$
\begin{figure}[h!]
\centering
\includegraphics[scale=.3]{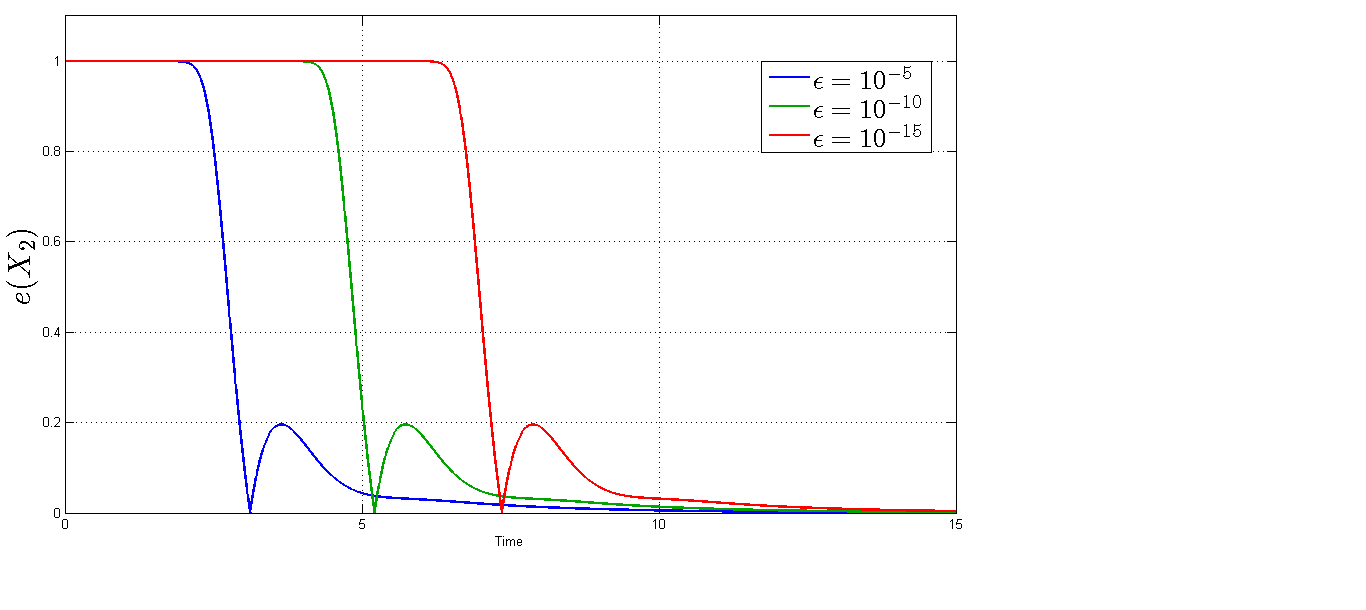}
\caption{Error between the current state and the desired reference for different initial conditions close to the undesired equilibria $\mathcal{R}_a(\pi,e_1)$-Smooth controller-}
\label{smooth_fig}
\end{figure}

\begin{figure}[h!]
\centering
\includegraphics[scale=.32]{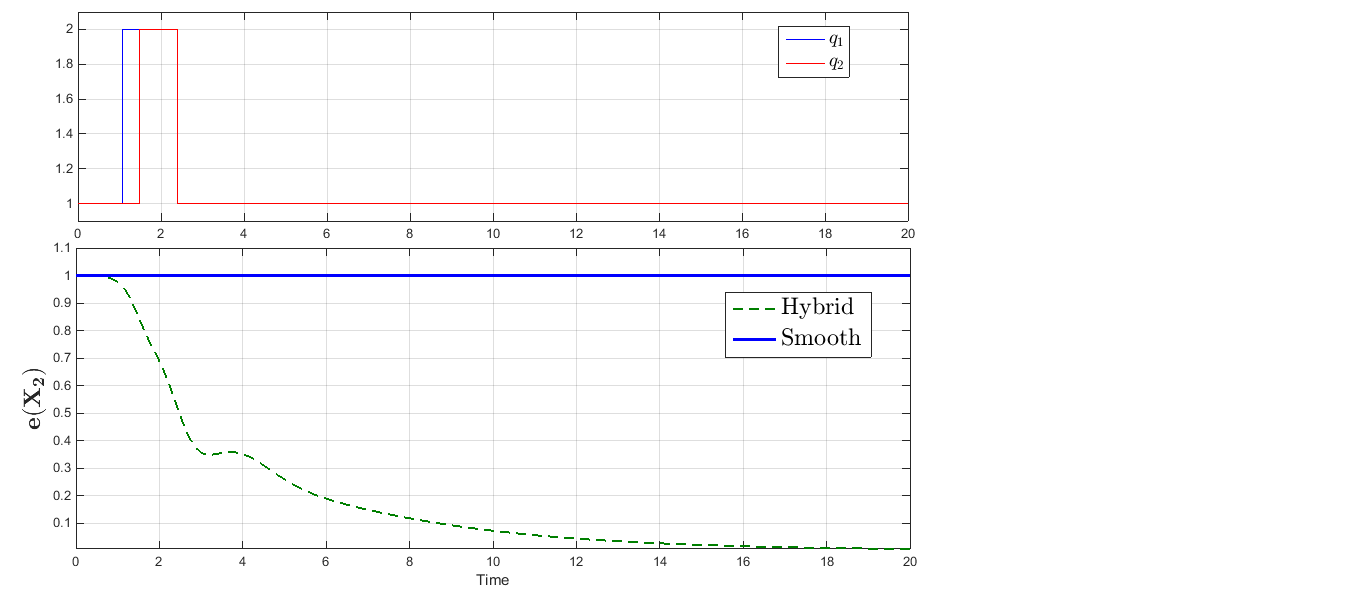}
\caption{Comparison between smooth and hybrid feedback responses with the initial condition $X_2(0)=\mathcal{R}_a(\pi,e_1)$}
\label{Smooth_hybrid_fig}
\end{figure}

To simulate the worst case for the smooth controller, the initial conditions for the rotational dynamics are taken as follows: $\omega(0)=[0,0,0]^{\top}$ and $R(0)=\mathcal{R}_a\left(\pi+\epsilon,\cos(\epsilon/2)e_1\right)$ for some $\epsilon\ll 1$. Thus, under the smooth feedback (\ref{smooth}), the closed-loop system starts sufficiently close to the undesired equilibria $X_1=X_2=\mathcal{R}(\pi,e_1)$. Figure \ref{smooth_fig} shows the evolution of the error $e(X_2)$ with respect to time.  At early times, the convergence is slower for closer initial conditions to the undesired equilibria (smaller choice of $\epsilon$). This phenomenon is the main drawback of the smooth controller. On the other hand, for $\epsilon=0$, Figure \ref{Smooth_hybrid_fig} depicts how the hybrid controller reacts immediately to correct its offset rotation, whereas the smooth controller does not react at all, being seemingly unable to correct its rotation.
\begin{figure}[h!]
\centering
\includegraphics[scale=.32]{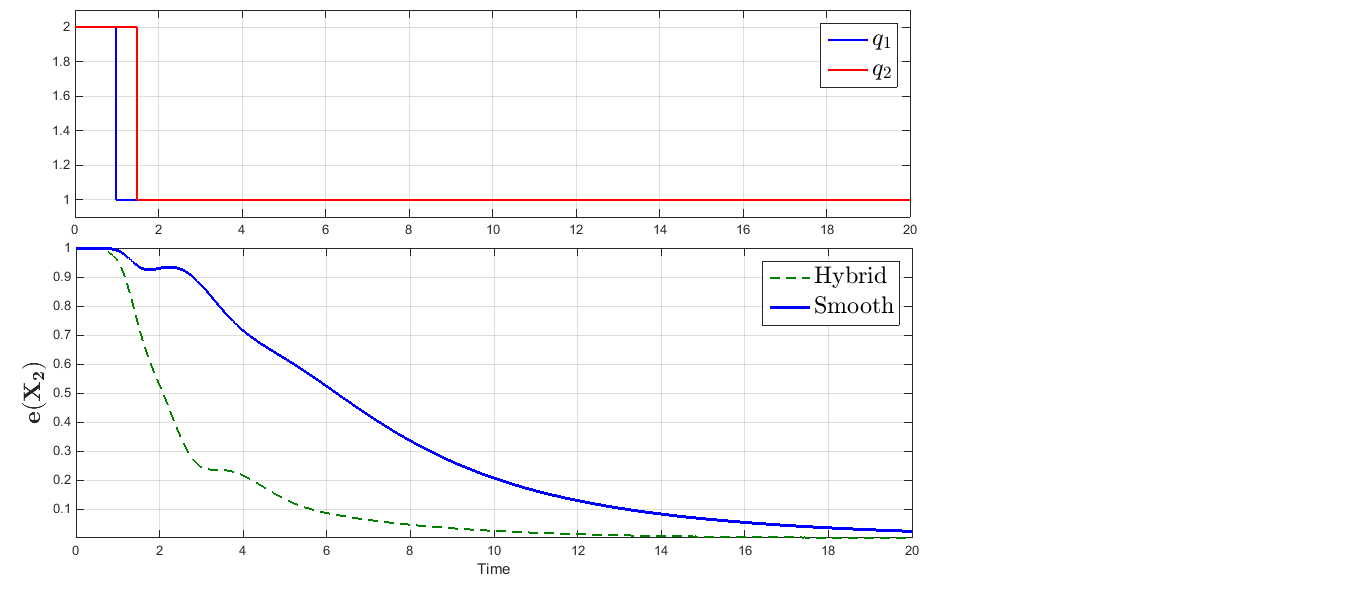}
\caption{Comparison between smooth and hybrid feedback responses with the initial condition $X_2(0)=\mathcal{R}_a(\pi,e_1)\mathcal{R}_a\left(\vartheta_1,u_1\right)^{\top}.$}
\label{Hybrid_fig}
\end{figure}
\begin{figure}[h!]
\centering
\includegraphics[scale=.29]{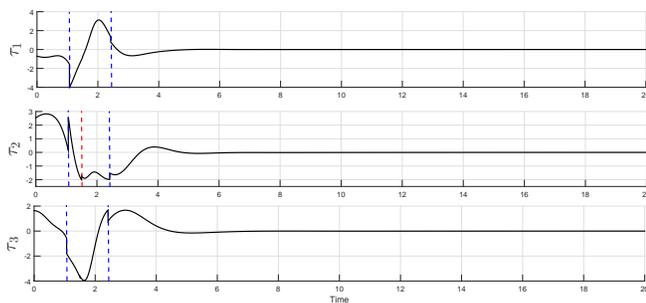}
\caption{Plot of the torque applied by the hybrid feedback with the initial condition $X_2(0)=\mathcal{R}_a(\pi,e_1)$.}
\label{tau1}
\end{figure}

For a second comparison, we changed the initial rotation matrix to
$$
R(0)=\mathcal{R}_a(\pi,e_1)\mathcal{R}_a\left(\vartheta_1,u_1\right)^{\top},
$$
where
$
\vartheta_1=2\arcsin\left(\frac{-1+\sqrt{1+16\lambda^{W_1}_1 (k_{11})^2\Delta(v_1,u_1)}}{4k_{11}\Delta(v_1,u_1)}\right)$, thus $\vartheta_1\simeq	0.47,
$
so as to start from one of the critical points of the potential function $\mathcal{U}_1(X_1,q_1)$ (see equation \eqref{critical_point2}).

Figure \ref{Hybrid_fig} depicts the performance of the proposed hybrid feedback law. Starting from an initial configuration $q(0)=(1,1)$, the system immediately jumps to the configuration $q=(2,2)$ since the initial condition $(X(0),q(0))$ lies inside the jump set $D$. It is shown in Figure \ref{Hybrid_fig} that the hybrid controller still achieves better performance than the continuous controller for this particular initial condition. 

In Figure \ref{tau1}, we give the plot of the torque (control input) applied by the hybrid controller in the first case of $X_2(0)=\mathcal{R}_a(\pi,e_1)$. 	We observe that the torque is ``quasi-smooth" with only three discontinuities which occur during the first few seconds of the control in order to avoid the critical points. The second jump (red) affects only the second component of the torque vector.

\section{Conclusion}\label{Conclusion}
Synergistic potential functions are instrumental in the design of hybrid control systems on $SO(3)$ that achieve global asymptotic stability results. This paper presented a systematic approach to generate synergistic potential functions on $SO(3)$ via angular warping. By introducing a new warping angle function, the synergistic gap-necessary for the implementation of the hybrid controller- was explicitly computed. The feasibility of the synergism conditions and the maximization of the synergy gap are discussed. We also proposed a hybrid attitude stabilization control scheme without velocity measurements relying only on inertial vector measurements. The proposed control scheme leads to global asymptotic stability results. We presented some simulation results that illustrate the advantage of the hybrid control scheme over the standard continuous feedback strategies. 
\bibliographystyle{IEEEtran}
\bibliography{IEEEabrv,Critical_points}
\end{document}